\documentclass[12pt]{article}

\title{Geodesic connectedness and conjugate points in GRW spacetimes.}

\author{Jos\'{e} Luis Flores and Miguel S\'anchez \thanks{Research partially
supported by a 
MEC grant number PB97-0784-C03-01.} \\ Depto. Geometr\'{\i}a y
Topolog\'{\i}a, Fac. Ciencias, Univ. Granada, \\ Avda. Fuentenueva s/n, 18071-Granada, Spain.}

\usepackage{graphicx}
\usepackage{amssymb}
\usepackage{latexsym}
\usepackage[latin1]{inputenc}

\newfont{\bb}{msbm10 at 12pt}
\newfont{\bt}{msbm10 at 8pt}
\newfont{\btt}{msbm10 at 5pt}
\newfont{\kt}{cmti10 at 6pt}
\newfont{\bg}{msbm10 at 16pt}

\newtheorem{lemma}{Lemma}
\newtheorem{remark}{Remark}
\newtheorem{theorem}{Theorem}
\newtheorem{proposition}{Proposition}
\newtheorem{corollary}{Corollary}
\newtheorem{convention}{Convention}

\newtheorem{definition}{Definition}

\font\ddpp=msbm10 scaled \magstep 1 
\def\R{\hbox{\ddpp R}}    
\def\N{\hbox{\ddpp N}}    

\def\be{\begin{equation}}
\def\ee{\end{equation}}

\begin{document}
\maketitle


\begin{abstract}

Given two points of a Generalized Robertson-Walker spacetime, the existence, multiplicity and causal character of geodesics connecting them is characterized. Conjugate points of such geodesics are related to conjugate points of geodesics on the fiber, and Morse-type relations are obtained. Applications to bidimensional spacetimes and to GRW spacetimes satisfying the timelike convergence condition are also found.

\end{abstract}

\newpage 
\section{Introduction}

Recently, geodesic connectedness of Lorentzian manifolds has been widely studied, and some related questions appear with great interest, among them: to determine the existence, multiplicity and causal character of geodesics connecting two points, to study their conjugate points and to find Morse-type relations. These questions have been answered, totally or partially, for stationary or splitting manifolds (see, for example, \cite{Ma}, \cite{BF}, \cite{GMP}, \cite{GMPT}, \cite{Uh}). Our purpose is to answer them totally in the class of  Generalized Robertson-Walker (GRW) spacetimes. 

GRW spacetimes (see Section \ref{s2} for precise definitions) are warped products $(I\times F, g^f=-dt^2 + f^2 g)$ which generalizes Robertson-Walker ones because no assumption on their fiber is done, and they have interesting properties from both, the mathematical and the physical point of view \cite{ARS}, \cite{Sa97}, \cite{Sa98}, \cite{Sa99}. GRW spacetimes are also particular cases of {\it multiwarped spacetimes}, whose geodesic connectedness has been recently studied by using a topological method  \cite{FS}. They can be also seen as splitting type manifolds, studied in \cite[Chapter 8]{Ma}, or as a type of Reissner Nordstr\"om Intermediate spacetimes, studied in \cite{Gi}, \cite{GM}. Nevertheless, we will see here that the results  for GRW spacetimes can be obtained in a simpler approach, and are sharper. In fact, we will develop the following direct point of view.

Given a geodesic  $\gamma(t)=(\tau(t), \gamma_F(t))$ of the GRW spacetime, the component $\gamma_F$ is a pregeodesic of its fiber. So, if $d\tau/dt$ does not vanish, we can consider the reparameterization $\gamma_F(\tau)$, and $\gamma$ will cross a point $z_0=(\tau_0,x_0)$ if and only if $x_0=\gamma_F(\tau_0)$. This simple fact yields a result on connectedness by timelike and causal geodesics \cite[Theorems 3.3, 3.7]{Sa97}. For spacelike geodesics, the reparameterization $\gamma_F(\tau)$ may fail. This problem can be skipped sometimes  by simple arguments on continuity  \cite[Theorem 3.2]{Sa98}, but we will study systematically it in order to solve completely the problem of geodesic connectedness. Moreover, this will be also the key to solve the other related problems (multiplicity,  conjugate points, etc.) 

After some preliminaries in Section \ref{s2}, we state the conditions for geodesic connectedness in Section \ref{s3}. In fact, we give three Conditions (A), (B), (C) of increasing generality, and a fourth Condition (R) which covers a residual case. All these conditions are imposed  on the warping function $f$; on the fiber, we assume just a weak condition on convexity (each two points $x_0, x'_0$ can be joined by a minimizing $F$--geodesic $\hat \gamma_F$), which is known to be completely natural (see  \cite[Remark 3.2]{Sa97}). 
These Conditions are somewhat cumbersome, because they yield not only  sufficient but also  {\em necessary} hypotheses for geodesic connectedness; however, they yield very simple sufficient conditions. For example, (Lemmas \ref{l00}, \ref{l3}) {\it if the GRW spacetime is {\em not} geodesically connected then $f$ must admit a limit at some extreme of the interval $I=(a,b)$; if this extreme is $b$ (resp. $a$) then $f'$ must be strictly positive (resp. negative) in a non-empty subinterval $(\bar b,b) \subseteq (a,b)$ (resp. $(a,\bar a) \subseteq (a,b)$)} (moreover, in this case Table 1 can be used). 
Condition (A) collects when the warping function $f$ has a ``good behavior'' at the extremes of $I=(a,b)$, in order  to obtain geodesic connectedness. This condition is equal to the one obtained in \cite{FS} for  multiwarped spacetimes; nevertheless, we will reprove it because a simpler proof is now available and the ideas in this proof will be used in the following more general conditions.  Condition (B) takes into account that, when the diameter of the fiber is finite, even a ``not so good'' behaviour of $f$ at a extreme, say $b$,  may allow  the following situation: a fixed point $z_0=(\tau_0,x_0)$ can be connected to  $z'_0=(\tau'_0,x'_0)$, where $\tau'_0$ is close enough to $b$, by means of a geodesic $\gamma(t)=(\tau(t), \gamma_F(t))$ such that $\tau(t)$ points out from $\tau_0$ to $b$ and, perhaps, ``bounces'' close to $b$. Condition (C) takes into account that, even when Condition (B) does not hold, the following situation in  the previous  case  may hold: a geodesic which points out from $\tau_0$ to the extreme $a$, bounces close to $a$ and comes back towards $b$, may connect $z_0$ and $z'_0$. Condition (C) is shown to be the more general condition for geodesic connectedness, except in a case: if the limit of $f$ at both extremes $a, b$ is equal to the supremum of $f$, and this supremum is not reached at $I$,   then $z_0$ and $z'_0$ perhaps could be joined by geodesics which bounces many times close to $a$ and $b$. 
Examples of the strict implications between the different Conditions are provided. In Section \ref{s3} we also state our results on existence of connecting geodesics, which are proven  in Section \ref{s4}: 
\begin{enumerate}
\item  either condition (C) or (R) is sufficient for geodesic connectedness (Theorem~\ref{t0})
\item if we assume a stronger condition of  convexity  on the fiber (each geodesic $\hat \gamma_F$ above is assumed to be the only geodesic which connects $x_0, x'_0$) then one of this two conditions (C) or (R) is also necessary (Theorem \ref{t1}); the necessity of this stronger convexity assumption is also discussed, 
\item under Condition (A) (or, even in some cases (B)), if the topology of $F$ is not trivial then each two points $z_0=(\tau_0,x_0), z'_0=(\tau'_0,x'_0)$
can be joined by infinitely many spacelike geodesics 
(Theorem \ref{t2}) and 
\item for causal geodesics: (i) if $z_0$ and $z'_0$ are causally related then there exist a causal geodesic joining them (this result was previously proven in \cite{Sa97}), (ii) if $z_0$ and $z'_0$  are not conjugate (or even if just  $x_0$ and $x'_0$ are not conjugate, which will be shown to be less restrictive), then there are at most finitely many timelike geodesics joining them (Theorem \ref{t2})  and (iii) if the fiber is strongly convex then there exist at most one connecting causal geodesic (Theorem \ref{t1}). 
\end{enumerate}

This machinery is used in Section \ref{s5} to obtain a precise relation between the conjugate points of a geodesic in $I\times F$ and its projection on $F$ (Theorem \ref{t-3}, Corollary \ref{cUh}). From this result, 
Morse-type relations, which relate the topology of the space of curves joining two non-conjugate points and the Morse indexes of the geodesics joining them, are obtained (see Corollary \ref{cM} and the discussion above it). We remark that Morse indexes are  defined here in the geometrical sense ``sum of the orders of conjugate points'' because, for any spacelike geodesic,  its index form is positive definite and negative definite on infinite dimensional subspaces (if dim$F>1$). About this kind of problem, the following previous references should be taken into account. Conjugate points of null geodesics in globally hyperbolic spacetimes were studied by Uhlenbeck \cite{Uh}, and we also make some remarks in Section \ref{s5} relating  our results. In a general setting, conjugate points on spacelike geodesics were studied by Helfer \cite{He}, who also considered the {\it Maslov index} of a geodesic. He showed that these conjugate points may have very different properties to conjugate points for Riemannian manifolds (unstability,  non-isolation...), but these problems can be skipped in our study. In  \cite{BM} (see also \cite[Section 5]{Ma}),  an attempt to obtain a Morse Theory for standard stationary manifolds is carried out, and in \cite{GMPT}, an index theorem (in terms of the Maslov index)  appliable in particular to stationary manifolds, is obtained. On the other hand,  some recent articles study a Morse theory for timelike or lightlike geodesics joining a point and a timelike curve, see \cite{GMPb} and references therein. Tipically, these results are stated for strongly causal spacetimes (including so all GRW spacetimes), and they need an assumption on {\it coercivity} which not necessarily holds under our hypotheses. It is not difficult to check that our results are also appliable to face this problem.

In Section \ref{s6} we particularize the previous results to two cases. First, Subsection \ref{s6.1}, when the fiber is also an interval of $\R$. In this bidimensional case, the opposite metric $-g^f$ is standard static, and we reobtain and extend the Theorem in \cite{BGM}. We recall that the proof in this reference is obtained by a completely different method, which relies in the function spectral flow on a geodesic (see the Remark (2) to Theorem \ref{t-3} for noteworthy comments about this approach). Finally, in Subsection \ref{s6.2} we consider the case Ric$(\partial_t , \partial_t )\geq 0$. This condition is natural from a physical point of view. In fact, the stronger condition Ric$(v,v)\geq 0$ for all timelike $v$, is called the timelike convergence condition, and says that gravity, on average, attracts. Condition Ric$(\partial_t , \partial_t )\geq 0$ is equivalent to $f''\leq 0$, and this inequality implies Condition (A) if $f$ cannot be continuously extended to positive values at any extreme. Corollary \ref{c00} summarizes our results in this case. We finish with an extension, in our ambient, of a result in \cite{Uh} (Corollary \ref{c6}).

\section{Preliminaries} \label{s2}

Let $ (F,g) $ be a Riemannian manifold, $ (I,-d\tau ^{2})$  an open interval of $\R{}$ 
with $I=(a,b)$ and its usual metric reversed, and $f>0$ a  smooth function on $I$. A {GRW spacetime} with 
base $(I,-d\tau^{2})$, fiber $(F, g)$ and warping function $f>0$ is the product 
manifold $I\times F$ endowed with the Lorentz metric:

\begin{equation}\label{e-1}g^{f}=-\pi_{I}^{*}d\tau^{2}+(f\circ\pi_{I})^{2}\pi^{*}_{F}g\equiv-d\tau^{2}+f^{2}g\end{equation}
where $\pi_{I}$ and $\pi_{F}$ are the natural projections of $I\times F$ onto I and $F$, respectively, and will be omitted when there is no possibility of confusion.

A Riemannian manifold will be called weakly convex if any two of its points  can be joined by a geodesic which minimize the distance; if, in addition, this geodesic is the only one which joins the two points it will be called strongly convex (recall that these names does {\it not} coincide with those in \cite{Sa97}). Of course, if the Riemannian manifold $(F,g)$ is complete then it is weakly convex by the Hopf-Rinow theorem, but the converse is not true; a detailed study of when a (incomplete) Riemannian manifold is weakly convex can be seen in \cite{BGS}. 
It is well-known that Cartan-Hadamard manifolds (i.e. complete, simply connected and with non-positive curvature) are strongly convex and, 
of course, so are locally all Riemannian manifolds; more results on strong convexity can be seen in \cite{GMP}.
We will denote by $d$ the distance on $F$ canonically associated to the Riemannian metric $g$,  and by $diam(F)$ its diameter (the supremum, possibly infinity, of the $d$-distances between points of $F$).  

Given a vector $X$ tangent to $I\times F$ we will say that $X$ is timelike (resp. lightlike, causal, spacelike) if $g^{f}(X,X)<0$ (resp. $=0, \leq 0, >0$); the timelike vector field $\partial /\partial\tau$ fixes the canonical future orientation in $I\times F$. Given $z,z'\in I\times F$, we will say that they are causally [resp. chronologically] related if they can be joined ($z$ with $z'$ or viceversa) by a future-pointing non-spacelike [resp. timelike]  piecewise smooth curve.

 Let $\gamma: {\cal J}\rightarrow I\times F$, $\gamma(t)=(\tau(t),\gamma_{F}(t))$ be a (smooth) curve on the interval ${\cal J}$. It is well-known that $\gamma$ is a geodesic with respect to $g^{f}$ if and only if 

\begin{equation} \label{e-2} \frac{d^{2}\tau}{dt^{2}}=-\frac{c}{f^{3}\circ\tau}\cdot\frac{df}{d\tau}\circ\tau \end{equation}

\begin{equation} \label{e-3} \frac{D}{dt}\frac{d\gamma_{F}}{dt}=-\frac{2}{f\circ\tau}\cdot\frac{d(f\circ\tau)}{dt}\cdot\frac{d\gamma_{F}}{dt} \end{equation}                                                                                                                                                                                                                                                                                                                                                    on ${\cal J}$, where $D/dt$ denotes the covariant derivate associated to $\gamma_{F}$ and $c$ is the constant $(f^{4}\circ\tau)\cdot g(d\gamma_{F}/dt,d\gamma_{F}/dt)$. 
From (\ref{e-2}), 
\begin{equation}\label{e-4}\frac{d\tau}{dt}=\epsilon\sqrt{-D+\frac{c}{f^{2}\circ\tau}}\end{equation}
with $D=g^{f}(d\gamma /dt,d\gamma /dt)$ and $\epsilon\in \{\pm1\}$.
Note that if $c=0$  then $d^{2}\tau /dt^{2}\equiv 0$, that is, the geodesics on the base  $I$ are naturally lifted to geodesics of  the GRW spacetime, as in any warped product. For all the other geodesics, it is natural to normalize choosing them with $c=1$. This normalization will be always chosen, except in Section \ref{s1} where the formulas will be explicitly taken with a different normalization. All geodesics will  be also assumed inextendible, that is, with a maximal domain.

By equation (\ref{e-3}), each (non-constant) $\gamma_{F}$ is a pregeodesic of $(F,g)$, so if we consider the reparametrization $\hat{\gamma}_{F}(r)=\gamma_{F}(t(r))$ where
                                                                             \begin{equation}\label{e-5}\frac{dt}{dr}=
f^{2}\circ\tau\circ t\end{equation}
(in a maximal domain) we obtain that $\hat{\gamma}_{F}$ is a geodesic of $(F,g)$ being

\begin{equation}\label{e-6}
g(\frac{d\hat{\gamma}_{F}}{dr},\frac{d\hat{\gamma}_{F}}{dr}) = 1.\end{equation}.

From now on, we will assume that $(F,g)$ is weakly convex for any result where geodesic connectedness is involved; such assumption has proven to be completely natural  \cite{Sa97}, \cite{Sa98}, \cite{FS}. In fact, as an immediate consequence of (\ref{e-2}) and (\ref{e-3}) we have:

\begin{lemma}\label{l-2}
There exists a geodesic joining 
$z_0=(\tau_{0},x_{0})$ and $z_0'=(\tau_{0}',x_{0}'), \tau_0 \leq \tau_0'$ if $\tau_{0} =\tau_{0}'$ and  $\frac{d}{d\tau}\frac{1}{f^{2}}(\tau_{0})=0$.
\end{lemma}

Now, the case when the geodesic $\hat{\gamma}_{F}$ can be reparameterized by using $\tau \in (a_*,b_*)$ as a parameter (for some   interval $(a_*,b_*)$) 
will be considered. Putting $\tilde \gamma (\tau) \equiv \hat \gamma  \circ r (\tau )$ we have $d\tilde \gamma / d\tau = h^\epsilon_D(\tau) \cdot d\hat \gamma / dr $ where $h^\epsilon_D \equiv h^\epsilon: (a_{\star},b_{\star})\subseteq I\rightarrow \R{}$     
 is defined as:
\be
\label{h}
h^{\epsilon}=\epsilon f^{-2}(-D+\frac{1}{f^{2}})^{-1/2},
\ee
and $\epsilon = \pm 1$. When this reparameterization can be done in such a way that if $\tau$ goes from $\tau_0 $ to $\tau'_0$  the integral of $h_\epsilon$ is exactly equal to the distance between $x_0, x'_0 \in F$, then a geodesic joining $(\tau_{0},x_{0})$ and $(\tau_{0}',x_{0}')$ can be constructed, yielding so:

\begin{lemma}
\label{l-1}
 There exists a geodesic connecting $z_0=(\tau_{0},x_{0})$ and $z_0'=(\tau_{0}',x_{0}'), \tau_0 \leq \tau_0'$ if there is a constant  $D\in \R{}$ ($D=g^{f}(\frac{d\gamma}{dt},\frac{d\gamma}{dt})$) such that,

(i) Either $\frac{1}{f^{2}(\tau_{0})}\neq D$ or if this equality holds then $\frac{d}{d\tau}\frac{1}{f^{2}}(\tau_{0})\neq 0$ and

(ii) the maximal domain $(a_{\star},b_{\star})$ of $h_\epsilon$ includes $(\tau_{0},\tau_{0}')$  and
\begin{equation}\label{e-8}\int_{\tau_{0}}^{\tau_{0}'}h^{\epsilon}=L .\end{equation}
where $L=d(x_{0},x'_{0})$.
 \end{lemma}
(In this case $\epsilon =1$.) When the reparameterization $\tilde \gamma (\tau)$ fails then the points where the denominator of $h^\epsilon$ goes to zero must be especially taken into account. Firstly, we will specify the maximal domain of $h^\epsilon$. Fix  $D\in \R{}$ such that $1/f^{2}(\tau_0)\geq D$ and consider the subsets
\begin{equation}\label{e-22}A_{+}=\{\tau\in (a,b): \tau_{0}\leq \tau, \frac{1}{f^{2}(\tau)}=D\}\cup\{b\}\end{equation}
\begin{equation}\label{e-23}A_{-}=\{\tau\in (a,b): \tau_{0}\geq \tau, \frac{1}{f^{2}(\tau)}=D\}\cup\{a\}\end{equation}
 Define $a_{\star}\equiv a_{\star}(D)$, $b_{\star}\equiv b_{\star}(D)$ by:

\begin{equation}\label{e-23.5} \hbox{If} \; \frac{d}{d\tau}\frac{1}{f^{2}}(\tau_{0})\left \{ \begin{array}{c} >0 \\ <0 \\ =0 \end{array} \right. \; \hbox{then} \; b_{\star}= \left \{ \begin{array}{c}  min(A_{+}-\{\tau_{0}\})  \\ min(A_{+}), \\ min(A_{+}) \end{array} \right.  a_{\star}= \left \{ \begin{array}{c} max(A_{-}) \\ max(A_{-}-\{\tau_{0}\}) \\ max(A_{-}) \end{array} \right. \end{equation}
  
Now, it is not difficult to check that  Lemma \ref{l-1} also holds if we assume the following convention for the integral (\ref{e-8}).

 \begin{convention}\label{cv1} {\rm From now on integral (\ref{e-8}) will be understood in the following generalized sense: for $\epsilon=1$, if $\int_{\tau_{0}}^{b_{\star}}h^{\epsilon=1}\geq L$, then the first member of (\ref{e-8}) denotes the usual integral and we will also follow the notation $_{+[0]}\int_{\tau_{0}}^{\tau'_{0}}h^{\epsilon}$; otherwise and if $b_{\star}\neq b$, we can follow integrating, by reversing the sense of integration (recall $\tau_0' \leq b_{\star}$) and, if $\int_{\tau_{0}}^{b_{\star}}h^{\epsilon=1}-\int_{b_{\star}}^{a_{\star}}h^{\epsilon=1}\geq L$  then the first member of (\ref{e-8}) means $\int_{\tau_{0}}^{b_{\star}}h^{\epsilon=1}-\int_{b_{\star}}^{\tau_{0}'}h^{\epsilon=1}$ which we denote by $_{+[1]}\int_{\tau_{0}}^{\tau'_{0}}h^{\epsilon}$. If this last inequality does not hold and $a_{\star}\neq a$, then the procedure must follow reversing the sense of integration ($\tau'_{0}\geq a_{*}$) as many times as necessary in the obvious way. Analogously, when $\epsilon=-1$, first member of (\ref{e-8}) means either $\int_{\tau_{0}}^{\tau_{0}'}h^{\epsilon=-1}\equiv _{-[0]}\int_{\tau_{0}}^{\tau'_{0}}h^{\epsilon}$ (in this case if  $\tau_{0} < \tau_{0}'$ the integral is negative; so  equality (\ref{e-8})  cannot hold) or  $\int_{\tau_{0}}^{a_{\star}}h^{\epsilon=-1}-\int_{a_{\star}}^{\tau_{0}'}h^{\epsilon=-1}\equiv _{-[1]}\int_{\tau_{0}}^{\tau'_{0}}h^{\epsilon}$  or $\int_{\tau_{0}}^{a_{\star}}h^{\epsilon=-1}-\int_{a_{\star}}^{b_{\star}}h^{\epsilon=-1}+\int_{b_{\star}}^{\tau_{0}'}h^{\epsilon=-1}\equiv _{-[2]}\int_{\tau_{0}}^{\tau'_{0}}h^{\epsilon}$, etc.}
\end{convention}

\begin{remark}\label{r0} {\rm From (\ref{h}), fixed $\epsilon\in\{\pm 1\}$, for each $D\in \R{}$ we have at most one $\tau'_{0}$ such that equation (\ref{e-8}) holds, possibly under Convention \ref{cv1}. Let us introduce the parameter $K\equiv K(D,\epsilon)$ by means of $K=\frac{1}{f^{2}(\tau_{0})}-D$ if $\epsilon =1$, $K=D-\frac{1}{f^{2}(\tau_{0})}$ if $\epsilon =-1$. So, fixed $L$, a function $\tau(K)=\tau'_{0}$ is defined for $K$ in certain domain ${\cal D}$ of $\R{}$.}
\end{remark}

\section{Conditions for geodesic connectedness} \label{s3}

Now, we are ready to stablish four conditions (Conditions (A), (B), (C), (R)) on the warping function $f$ which, independently, ensure the geodesic connectedness of the GRW spacetime (Lemma \ref{l2}, Lemma \ref{l3} and Theorem \ref{t0}). Roughly, Condition (A) implies not only the geodesic connectedness but also that every $(\tau_{0},x_{0})\in I\times F$ can be joined with any point $(\tau'_{0},x'_{0})$ with 
$\tau'_0$ close enough to $b$ (resp. $a$) by means of a geodesic 
 $(\tau(t), \gamma_{F}(t))$ with $d\tau /dt>0$ (resp. $<0$) near $\tau'_{0}$. 
Condition (B) is weaker than Condition (A), and implies not only geodesic connectedness but also that if Condition (A) does not hold at $b$ (resp. $a$) then any  $(\tau_{0},x_{0})\in I\times F$ can be joined with  a point $(\tau'_{0},x'_{0})$ with 
$\tau'_0$ close enough to $b$ (resp. $a$)
   by means of a geodesic with  $\epsilon=1$ (resp. $\epsilon=-1$), and perhaps using Convention \ref{cv1} once close to $\tau'_0$. Condition (C) is the most general condition for geodesic connectedness, which just drops a residual case covered by  Condition (R).

\begin{definition}\label{d-1} Let $f:(a,b)\rightarrow \R{}$ be a smooth function and let $m_{b}=lim\; inf_{\tau\rightarrow b}f(\tau)$ [resp. $m_{a}=lim\; inf_{\tau\rightarrow a}f(\tau)$]. The extreme $b$ [resp. $a$] is a (strict) relative minimum of $f$ if:

(a) When $b<\infty$ [resp. $a>-\infty$], there exists $\epsilon>0$ such that if $0<\epsilon'<\epsilon$, then $f(b-\epsilon')>m_{b}$ [resp. $f(a+\epsilon')>m_{a}$].

(b) When $b=\infty$ [resp. $a=-\infty$], there exist $M>0$ such that if $M'>M$ then $f(M')>m_{b}$ [resp. $f(-M')>m_{a}$].
\end{definition}

{\bf Condition (A) for $f$}. Either  $1/f^{2}$ does not reach at $b$ [resp. $a$] a relative minimum in the sense of Definition \ref{d-1} or, otherwise,  $$\int_{c}^{b}f^{-2}(\frac{1}{f^{2}}-m_{b})^{-1/2}=\infty \quad [\mbox{resp.} \int_{a}^{c}f^{-2}(\frac{1}{f^{2}}-m_{a})^{-1/2}=\infty]$$ for some $c\in (a,b)$ {\em close to} $b$ (resp. $a$), i.e. $c\in (b-\epsilon,b)$ [resp. $c\in (a,a+\epsilon)$] if the extreme $b$ (resp. $a$) is finite or $c>M$ [resp. $c<-M$] if this extreme is infinite, where $\epsilon$ and $M$ are given in Definition \ref{d-1}.

The following definition is needed to state Condition (B). Recall that this definition is appliable just when Condition (A) does not hold.

\begin{definition}\label{d1}  Assume that the function $1/f^{2}$ reaches at $b$ [resp. at $a$] a relative minimum such that   $\int_{c}^{b}f^{-2}(\frac{1}{f^{2}}-m_{b})^{-1/2}<\infty$  [resp. $\int_{a}^{c}f^{-2}(\frac{1}{f^{2}}-m_{a})^{-1/2}<\infty$] for some $c\in (a,b)$. Then we define $$d_{b}=lim\; sup_{D\rightarrow m_{b}} \left(\int_{c}^{b_{\star}}f^{-2}(\frac{1}{f^{2}}-D)^{-1/2}\right)-\int_{c}^{b}f^{-2}(\frac{1}{f^{2}}-m_{b})^{-1/2}$$  $$[\mbox{resp.}\quad d_{a}=lim\; sup_{D\rightarrow m_{a}} \left( \int_{a_{\star}}^{c}f^{-2}(\frac{1}{f^{2}}-D)^{-1/2}\right)-\int_{a}^{c}f^{-2}(\frac{1}{f^{2}}-m_{a})^{-1/2} ]$$  where $b^* \equiv b^*[D]$ (resp. $a^* \equiv a^*[D]$) is given by (\ref{e-23.5}).
\end{definition}

{\it Remark.} (1) Note that the uniform convergence  of $f^{-2}(\frac{1}{f^{2}}-D)^{-1/2}$ on compact subsets of $(a,b)$ when $D$ varies,  ensures that $d_{b}$ and $d_{a}$ are independent of $c$.

(2) It is easy to check that $d_{b}, d_{a}\geq 0$. 
As when $D\rightarrow m_b$ then $b^* \rightarrow b$, where $b$ is a relative minimum, it is clear that if there were continuity of the integrals with $D$ at $m_b$ [resp. $m_a$] then $d_b=0$ [resp. $d_a=0$]. But as we will see in the example below,  there exist cases in which the inequalities are strict, and $d_b, d_a$ can reach even the value $\infty$.

\vspace{3mm}

{\bf Condition (B) for $f$}. Either $1/f^{2}$ does not reach at $b$ [resp. $a$] a relative minimum or, otherwise, it verifies   either $\int_{c}^{b}f^{-2}(\frac{1}{f^{2}}-m_{b})^{-1/2}=\infty$  for some $c\in (a,b)$ as in Condition (A), or $2d_{b}\geq diam(F)\in (0,\infty]$  
(resp. either  $\int_{a}^{c}f^{-2}(\frac{1}{f^{2}}-m_{a})^{-1/2}=\infty$ or $2d_{a}\geq diam(F)\in (0,\infty]$).

\vspace{2mm}
Obviously Condition (A) implies Condition (B), but  the converse is not true as the following example shows. 

{\it Example.} Consider the function $1/g^{2}(\tau) =1-\tau$ defined on $(0,1)$. Modify this function smoothly on $\{I_{n}\}_{n\in \N{}}, \, I_{n}=(a_{n},b_{n}), \, a_{n},b_{n}\rightarrow 1, \, a_n<b_n<a_{n+1}$ in such a way that the modified function $1/f^{2}$ satisfies $1/f^{2}>1/g^{2}$ on $I_{n}, \, \forall n\in \N{}$ and   
\begin{equation} \label{eee} 
\int_{0}^{b_{n*}}f^{-2}(\frac{1}{f^{2}}-D_{n})^{-1/2}\geq 2L
\end{equation}  
where  $\int_{0}^{1}f^{-1}=L$ and $D_n$ is chosen decreasing to 0 and such that $b_{n*}=\frac{a_{n}+b_{n}}{2}$; this is possible by taking $1/f^{2}$ with derivative small enough in $(a_{n},\frac{a_{n}+b_{n}}{2})$ (for example, if this derivative vanishes at $\frac{a_{n}+b_{n}}{2}$ the integral (\ref{eee}) will be infinite). Then, as $d_{b}\geq L$, it is sufficient to take $(F,g)$ such that  $2d_{b}\geq diam(F)$ (see Fig. 1). 

\begin{lemma}\label{l00} 
If Condition (B) does not hold at $b$ (resp. $a$) then there exist $lim_{\tau\rightarrow b}f\in (0,\infty]$ (resp. $lim_{\tau\rightarrow a}f\in(0,\infty]$) and $f'>0$ on $(b-\delta,b)$ or $(M,\infty)$ [resp. $f'<0$ on $(a,a+\delta)$ or $(-\infty,-M)$] for some $\delta>0$ small or $M>0$ big.
\end{lemma}

{\it Proof.} Reasoning for $b<\infty$, assume that Condition (B) does not hold at $b$. Then $\frac{1}{f^{2}}$ reaches a relative minimum at $b$ and  $\int_{\tau_{0}}^{b}f^{-2}(\frac{1}{f^{2}}-m_{b})^{-1/2}<\infty$ for certain $\tau_{0}\in I$, see Definition \ref{d-1}. It is sufficient to prove that $f'>0$ on $(b-\delta,b)$. Otherwise, there exist a sequence $\{\overline{\tau}_{n}\}_{n\in \N{}}$, $\tau_{0}< \overline{\tau}_{n}\in I$, $\overline{\tau}_{n}\rightarrow b$ such that $f'(\overline{\tau}_{n})\leq 0$. If we choose a maximum $\tau_{n}$ of $f$ on  $[\tau_{0},\overline{\tau_{n}}]$,  then $f'(\tau_{n})=0$ for $n$ big enough. Thus, $\int_{\tau_{0}}^{b^{*}}f^{-2}(\frac{1}{f^{2}}-D_{n})^{-1/2}=\infty$ for $D_{n}=\frac{1}{f^{2}(\tau_{n})}$. The choice of $\tau_{n}$ implies that $D_{n}\rightarrow m_{b}$, which contradicts that $d_{b}<\infty$. $\Box$

{\it Remark.} If Condition (B) does not hold at $b$ (resp. $a$) then $lim_{\tau\rightarrow b}\frac{1}{f^{2}}=m_{b}$ (resp.  $lim_{\tau\rightarrow a}\frac{1}{f^{2}}=m_{a}$).

From  Lemma \ref{l00} it is natural to construct  Table 1, where it is assumed that $f$ is continuously extendible to $b$ (the Table for  $a$ would be analogous, but reversing the sign of the corresponding $\beta$).

The following definition, necessary to state Condition (C), is appliable when Condition (B) does not hold.

\begin{definition}\label{d2} Assume that  the function $1/f^{2}$ reaches at $b$ [resp. $a$] a relative minimum such that $2d_{b}<diam(F)$ and $m<m_{b}$ [resp. $2d_{a}<diam(F)$ and $m<m_{a}$] where $m$ is the infimum value of  $1/f^{2}$ in $(a,b)$. Choose $\tau_{0}\in (b-\epsilon,b)$ or $\tau_{0}>M$ [resp. $\tau_{0}\in (a,a+\epsilon)$ or $\tau_{0}<-M$] where $\epsilon, M$ are given in Definition \ref{d-1}. Then we define  
 $$i_{b}=inf_{D\in(m,m_{b}]}\{\int_{a_{\star}}^{b}f^{-2}(\frac{1}{f^{2}}-D)^{-1/2}\}$$  $$[\mbox{resp.}\quad i_{a}=inf_{D\in(m,m_{a}]}\{\int_{a}^{b_{\star}}f^{-2}(\frac{1}{f^{2}}-D)^{-1/2}\}]$$
where $b^* \equiv b^*[D]$ (resp. $a^* \equiv a^*[D]$) is given by (\ref{e-23.5}).

\end{definition}

Note that this definition is independent of the choice of $\tau_{0}$.

{\bf Condition (C)}. Either $1/f^{2}$ does not reach at $b$ [resp. $a$] a relative minimum or, otherwise,  either $\int_{c}^{b}f^{-2}(\frac{1}{f^{2}}-m_{b})^{-1/2}=\infty$ for some $c\in (a,b)$ as in Condition (A), or $2d_{b}\geq diam(F)$, or $d_{b}\geq i_{b}$ [resp. either  $\int_{a}^{c}f^{-2}(\frac{1}{f^{2}}-m_{a})^{-1/2}=\infty$  or $2d_{a}\geq diam(F)$ or $d_{a}\geq i_{a}$].

\vspace{3mm}

Again Condition (B) implies obviously Condition (C), and  a counterexample to  the converse is shown. 

{\it Example.} Let $1/f^{2}$ be the function in previous example.  We have that the smooth function $\bar f$ defined on $(-1/N,1)$ such that $\lim_{\tau\rightarrow -1/N} 1/\bar f^2 =0$,   $1/\bar f^2(0)=N+1$  and $1/\bar f^2(\tau) = N + 1/f^{2}(\tau) $ for  $\tau \in (0,1)$ satisfies that $i_{b}\leq d_{b}$ for $N$ big enough. Then it is sufficient to take $(F,g)$ such that $2d_{b}<diam(F)$ (see Fig. 2).
\vspace{3mm}

For the remaining residual case, we need the following definition, where Convention \ref{cv1} is explicitly used.

\begin{definition}\label{d2.5} Assume $\frac{1}{f^{2}}>m$ for $\tau\in (a,b)$ and $m_{a}=m_{b}=m$. Then  we define

\begin{equation}
\begin{array}{l}r_{i}^{n}(\tau_{0})=lim_{\epsilon\searrow 0}lim\; inf_{D\searrow m}\{_{(-)^{n}[n-1]}\int^{a_{*}}_{\tau_{0}}f^{-2}(\frac{1}{f^{2}}-D)^{-1/2}+\int_{a_{*}}^{b-\epsilon}f^{-2}(\frac{1}{f^{2}}-D)^{-1/2}\} \\ r_{s}^{n}(\tau_{0})=lim_{\epsilon\searrow 0}lim\; sup_{D\searrow m}\{_{(-)^{n}[n]}\int_{\tau_{0}}^{b_{*}}f^{-2}(\frac{1}{f^{2}}-D)^{-1/2}+\int_{b-\epsilon}^{b_{*}}f^{-2}(\frac{1}{f^{2}}-D)^{-1/2}\}  \\ l_{i}^{n}(\tau_{0})=lim_{\epsilon\searrow 0}lim\; inf_{D\searrow m}\{_{(-)^{n-1}[n-1]}\int_{\tau_{0}}^{b_{*}}f^{-2}(\frac{1}{f^{2}}-D)^{-1/2}+\int_{a+\epsilon}^{b_{*}}f^{-2}(\frac{1}{f^{2}}-D)^{-1/2}\} \\ l_{s}^{n}(\tau_{0})=lim_{\epsilon\searrow 0}lim\; sup_{D\searrow m}\{_{(-)^{n-1}[n]}\int^{a_{*}}_{\tau_{0}}f^{-2}(\frac{1}{f^{2}}-D)^{-1/2}+\int_{a_{*}}^{a+\epsilon}f^{-2}(\frac{1}{f^{2}}-D)^{-1/2}\}\end{array}
\end{equation}

for $n\geq 1$, and

\begin{equation}
\begin{array}{l}r_{i}^{0}(\tau_{0})=lim_{\epsilon\searrow 0}lim\; inf_{D\searrow m}\{\int_{\tau_{0}}^{b-\epsilon}f^{-2}(\frac{1}{f^{2}}-D)^{-1/2}\} \\ r_{s}^{0}(\tau_{0})=lim_{\epsilon\searrow 0}lim\; sup_{D\searrow m}\{\int_{\tau_{0}}^{b_{*}}f^{-2}(\frac{1}{f^{2}}-D)^{-1/2}+\int_{b-\epsilon}^{b_{*}}f^{-2}(\frac{1}{f^{2}}-D)^{-1/2}\}  \\ l_{i}^{0}(\tau_{0})=lim_{\epsilon\searrow 0}lim\; inf_{D\searrow m}\{\int_{a+\epsilon}^{\tau_{0}}f^{-2}(\frac{1}{f^{2}}-D)^{-1/2}\} \\ l_{s}^{0}(\tau_{0})=lim_{\epsilon\searrow 0}lim\; sup_{D\searrow m}\{\int_{a_{*}}^{\tau_{0}}f^{-2}(\frac{1}{f^{2}}-D)^{-1/2}+\int_{a_{*}}^{a+\epsilon}f^{-2}(\frac{1}{f^{2}}-D)^{-1/2}\}\end{array}
\end{equation}

\end{definition}

If some extreme of $I$ is infinite, previous definition must be understood in the natural way (see comments above formula (\ref{m8}))

Recall that $r_{i}^{0}(\tau_{0})=\int^{b}_{\tau_{0}}f^{-2}(\frac{1}{f^{2}}-m)^{-1/2}$ (resp. $l_{i}^{0}(\tau_{0})=\int_{a}^{\tau_{0}}f^{-2}(\frac{1}{f^{2}}-m)^{-1/2}$). It is also clear that $r_{i}^{n}(\tau_{0})\leq r_{s}^{n}(\tau_{0})$ (resp. $l_{i}^{n}(\tau_{0})\leq l_{s}^{n}(\tau_{0})$) and the sequence $\{r_{i}^{n}(\tau_{0})\}_{n\in \N{}}$ is strictly increasing to $\infty$ (resp. replacing $i$ or $r$ by $s$ or $l$).

{\bf Condition (R).} Assume $\frac{1}{f^{2}}>m$ for $\tau\in (a,b)$ and $m_{a}=m_{b}=m$, then 
$$[r_{i}^{0}(\tau_{0}),diam(F)]\subseteq \cup_{n\geq 0}[r_{i}^{n}(\tau_{0}),r_{s}^{n}(\tau_{0})] \; {\it and} \;      
 [l_{i}^{0}(\tau_{0}),diam(F)]\subseteq \cup_{n\geq 0}[l_{i}^{n}(\tau_{0}),l_{s}^{n}(\tau_{0})]$$
 for every $\tau_{0}\in I$.

\begin{remark}\label{r1}{\rm When $\frac{1}{f^{2}}>m$ for $\tau \in (a,b)$ and $m_{a}=m_{b}=m$, it is clear that Condition (C) holds if and only if Condition (B) holds; moreover Condition (R) is less restrictive than Condition (B). In fact, when Condition (A) holds then $r_{i}^{0}(\tau_{0})=\infty=l_{i}^{0}(\tau_{0})$ for all $\tau_{0}\in I$, thus Condition (R) is automatically satisfied. When Condition (A) does not hold then if $2d_{b}\geq diam(F)$ (i.e. Condition (B) holds at $b$) then $r_{s}^{0}(\tau_{0})\geq diam(F)$ for all $\tau_{0}\in I$ (and, thus Condition (R) holds).}         
\end{remark}

Condition (C) and Condition (R) provides us  accurate sufficient hypotheses for geodesic connectedness, as the following two theorems show. (For the sake of completeness, we also state  the result on connection by causal geodesics, already contained in \cite[Theorems 3.3, 3.7]{Sa97}).

\begin{theorem}\label{t0} 
Let $(I\times F,g^f=-d\tau^{2}+f^{2}g)$ be a GRW spacetime with weakly  convex fiber $(F,g)$. Then:

(i) Two points $z_0=(\tau_0, x_0), z'_0=(\tau'_0, x'_0), \tau_0 < \tau'_0$ are chronologically (resp. causally) related if and only if
$\int_{\tau_0}^{\tau'_0}f^{-1} > d_F(x_0,x'_0)$ (resp. $\geq d_F(x_0,x'_0)$) and, in this case, they 
 can be joined with at least one timelike (resp. non-spacelike) geodesic.

(ii) If Condition (C) or Condition (R) holds then the GRW spacetime is geodesically connected.
\end{theorem}

When the fiber is strongly convex, Condition (C) or Condition (R) becomes also necessary:

\begin{theorem}\label{t1} Let $(I\times F,g^f=-d\tau^{2}+f^{2}g)$ be a GRW spacetime with strongly convex fiber $(F,g)$. Then:

(i) Each two causally related points can be joined with exactly one (necessarily non-spacelike) geodesic. 

(ii) The GRW spacetime is geodesically connected if and only if either Condition (C) or Condition (R) holds.
\end{theorem} 

From its proof, it is clear the naturality of the strong convexity assumption. However we discuss, below the proof of Theorem \ref{t1}, what happens if just weak convexity is assumed. 

As a consequence of our technique, we also obtain the following result on multiplicity:

\begin{theorem}\label{t2} Let $(I\times F,g^f=-d\tau^{2}+f^{2}g)$ be a GRW spacetime with weakly convex fiber $(F,g)$ and assume that either Condition (A) or Condition (B) with $d_{a}$, $d_{b}$ (if defined) equal to infinity, holds.

Then there exist a natural surjective map between geodesics connecting $z_{0}=(\tau_{0},x_{0}),\; z'_{0}=(\tau'_{0},x'_{0})\in I\times F$ and F-geodesics connecting $x_{0}$ and $x'_{0}$.

Moreover, if $(F,g)$ is complete and $F$ is not contractible in itself then any $z_{0}, z'_{0}\in I\times F$ can be joined by means of infinitely many spacelike geodesics. If the corresponding $x_{0}, x'_{0}$ are not conjugate in $(F,g)$, then there are at most finitely many causal geodesics connecting $z_{0}, z'_{0}$ in $I\times F$.
\end{theorem}

{\it Remark.} From results in Section \ref{s1}, it will be clear that  to impose the non--conjugacy of $x_{0}, x'_{0}$ as above, is less restrictive than to impose the non--conjugacy of  $z_{0}, z'_{0}$. On the other hand, the completeness of the fiber in Theorem \ref{t2} can be replaced for a convexity assumption of the Cauchy boundary, as in \cite{BGS}.

\section{Proof of Theorems} \label{s4}

Consider a GRW spacetime  $(I\times F,-d\tau^{2}+f^{2}g)$  with weakly convex fiber $(F,g)$.
Fixed $\tau_0 \in I$ put 
\be \label{e-m}
m_{r} = Inf\{ 1/f^{2}(\tau) \mid \tau \in [\tau_{0},b)\}, \quad
m_{l} = Inf\{ 1/f^{2}(\tau) \mid \tau \in (a,\tau_{0}]\}.
\ee

\begin{lemma}\label{l1} 
Using the notation (\ref{e-23.5}), the function in $D$  $$\int_{\tau_{0}}^{b_{\star}}f^{-2}(\frac{1}{f^{2}}-D)^{-1/2}, b_*\equiv b_*(D) \quad  [\mbox{resp.}  \int_{a_{*}}^{\tau_{0}}f^{-2}(\frac{1}{f^{2}}-D)^{-1/2}, a_* \equiv a_*(D) ]$$ with values in $(0,\infty]$  is continuous when $D$ varies in  $(m_{r},\frac{1}{f^{2}(\tau_{0})})$  [resp.  $(m_{l},\frac{1}{f^{2}(\tau_{0})})$].
\end{lemma}

{\it Proof.}  We will check that every convergent sequence $\{D^{k}\}_{k\in \N{}}$, $D^{k}\rightarrow D^{\infty}$, $D^{\infty}\in (m_{r},\frac{1}{f^{2}(\tau_{0})})$  satisfies $\int_{\tau_{0}}^{b_{*}^{k}}f^{-2}(\frac{1}{f^{2}}-D^{k})^{-1/2}\rightarrow \int_{\tau_{0}}^{b_{*}^{\infty}}f^{-2}(\frac{1}{f^{2}}-D^{\infty})^{-1/2}$  (the case with $a_*$ is analogous). We can consider the following possibilities:

(i) If $\frac{d}{d\tau}\frac{1}{f^{2}}\mid_{b_{*}^{\infty}}\neq 0$ then the sequence of intervals $[\tau_{0},b_{*}^{k})$     converges to $[\tau_{0},b_{*}^{\infty})$ and the integrands converge uniformly on $[\tau_0, b_*^{\infty}-\delta]$ for $\delta>0 $ small, which implies the  convergence of the integrals in $[\tau_0, b_*^{\infty}-\delta]$. Thus, the result follows because the integrals on $[b_*^{\infty}-\delta, b_*^{\infty}]$ goes to zero when $\delta \rightarrow 0$.

(ii) If  $\frac{d}{d\tau}\frac{1}{f^{2}}\mid_{b_{*}^{\infty}}=0$  then the uniform convergence of $f^{-2}(\frac{1}{f^{2}}-D^{k})^{-1/2}$ to  $f^{-2}(\frac{1}{f^{2}}-D^{\infty})^{-1/2}$  on compact subsets of $[\tau_{0},b_{*}^{\infty})$  implies that  $\int_{\tau_{0}}^{b_{*}^{k}}f^{-2}(\frac{1}{f^{2}}-D^{k})^{-1/2}\rightarrow \infty=\int_{\tau_{0}}^{b_{*}^{\infty}}f^{-2}(\frac{1}{f^{2}}-D^{\infty})^{-1/2}$ . $\Box$

\vspace{2mm}
Recall that the integrals not necessarily varies continuously when $D= m_r , m_l$.

In what follows we will use the function $\tau(K)$ defined in Remark \ref{r0}, and follow the notation: 
$\tau^{-} = \tau(K^{-})$,  $\tau^{+} = \tau(K^{+})$

\begin{lemma}\label{l1.5} 
Consider $(\tau_{0},x_{0})\in I\times F$ and $x'_{0}\in F$ such that $d(x_{0},x'_{0})=L>0$. The function $\tau(K)$  is continuous on its domain ${\cal D}$. Moreover, if $\frac{d}{d\tau}\mid_{\tau=\tau_{0}}\frac{1}{f^{2}(\tau)}=0$ then $\tau(K)$ can be continuously extended to $K=0$ by $\tau(0)=\tau_0$. 

As a consequence, if  $[K^{-},K^{+}] \subset {\cal D}$ then we can connect $(\tau_{0},x_{0})$ with $[\tau^{-},\tau^{+}]\times \{x'_{0}\}$ (or $[\tau^{+},\tau^{-}]\times \{x'_{0}\}$).
\end{lemma}

{\it Proof.} Firstly, we will check that every convergent sequence $\{K^{n}\}_{n\in \N{}}$, $K^{n}\rightarrow K^{\infty}>0$ ($<0$ analogous), $K^{n},K^{\infty}\in {\cal D}$ for all $n$, satisfies that $\tau^{n}\rightarrow \tau^{\infty}$, where $\tau^{n} = \tau(K^{n})$, $\tau^{\infty} = \tau(K^{\infty})$. Assume first $K^{\infty} \neq 0$, then:

(i) If $\frac{d}{d\tau}\frac{1}{f^{2}}\mid_{b_{*}^{\infty},a_{*}^{\infty}}\neq 0$, then easily $a_{*}^{n}\rightarrow a_{*}^{\infty}$, $b_{*}^{n}\rightarrow b_{*}^{\infty}$ , so the proof follows from Lemma \ref{l1}.

(ii) If $\frac{d}{d\tau}\frac{1}{f^{2}}\mid_{b_{*}^{\infty}}=0$ then, as  $\int_{\tau_{0}}^{b_{*}^{\infty}}f^{-2}(\frac{1}{f^{2}}-D^{\infty})^{-1/2}=\infty$, we have $\tau^{\infty} < b_{\star}^{\infty}$ and the uniform convergence of the integrand on a compact set $[\tau_{0}, \tau^{\infty}+\delta]$ ($\delta>0$ small) proves the result.

(iii) If $b_{*}^{\infty}=b$  then again $\tau^{\infty}<b$ and the result follows from the convergence on  $[\tau_{0}, \tau^{\infty}+\delta]$.

 (iv) The remaining cases follows from combinations of the previous ones.

Now, consider the case that 
$ K^{\infty}=0 \in {\cal D}$
and (necessarily) $\frac{d}{d\tau}\frac{1}{f^{2}}\mid_{\tau_{0}}\neq 0$. Then it is easy to check that Lemma \ref{l1} 
can be extended to $D=\frac{1}{f^{2}(\tau_{0})}$, which implies the continuity of $\tau$ at 0. 

So, we have just to prove that if $\frac{d}{d\tau}\frac{1}{f^{2}}\mid_{\tau_{0}} = 0$,
then $\tau(K)$ can be continuously extended as $\tau(0)= \tau_0$. 
Fixed $\epsilon >0$, the limit of  $\int_{\tau_{0}}^{\tau_{0}+\epsilon}f^{-2}(\frac{1}{f^{2}}-D)^{-1/2}$ and  $\int_{\tau_{0}-\epsilon}^{\tau_{0}}f^{-2}(\frac{1}{f^{2}}-D)^{-1/2}$ (for the values of $D$ where they are well defined) are $\infty$ when $D\nearrow \frac{1}{f^{2}(\tau_{0})}$ (and, thus, $K \rightarrow 0$), from which the result follows. $\Box$

\begin{lemma}\label{l1.75} If  $K^{+}>0$ [resp. $K^{-}<0$] belongs to the domain ${\cal D}$ of $\tau(K)$ but  $K^{+}-\epsilon \geq 0$ [resp. $K^{-}+\epsilon \leq 0$] for some  $\epsilon >0$, does not belong, then we can connect $(\tau_{0},x_{0})$ with $(a,\tau^{+}]\times \{x'_{0}\}$ [resp. $[\tau^{-},b)\times \{x'_{0}\}$] 
by means of geodesics with $K \in (K^+-\epsilon, K^+]$ [resp. $K \in [K^-, K^- +\epsilon)$]. 

\end{lemma}

{\it Proof.} Reasoning for $K^{+}$, define $K_{0}=inf\{ K\leq K^{+}:
[K,K^{+}]\subseteq {\cal D}\}$. As $0 \leq K_{0}<K^{+}$, the fact that $K_{0}$ is the infimum implies that $b_{*}(D_{0})\neq b$ where $D_{0}=\frac{1}{f^{2}(\tau_{0})}-K_{0}$. Therefore, $\lim_{K\searrow K_0} \tau(K)=a$ (otherwise, it would contradict that $K_{0}$ is the infimum again) and the the result follows from the first assertion in Lemma \ref{l1.5}. $\Box$

\begin{lemma}\label{l.nuevo}
If the domain ${\cal D}$ contains $K^+ > 0$ and $K^- <0$, and
the inequality $\tau^{-}< \tau^{+}$ holds, then we can connect $(\tau_{0},x_{0})$ with, at least, $[\tau^{-},\tau^{+}]\times \{x'_{0}\}$  by choosing $K\in [K^{-},K^{+}]$.
\end{lemma}

{\it Proof.} If $\tau$ is defined in $[K^{-},K^{+}]$ then Lemma \ref{l1.5} can be applied. Otherwise, let $K_{0}\in (K^{-},K^{+})$ be such that $K_{0} \not\in {\cal D}$. If, say , $K_{0}\geq 0$  Lemma \ref{l1.75} can be applied to $K^{+}$. $\Box $

Now, a first result on geodesic connectedness can be stated.

\begin{lemma}\label{l2} A GRW spacetime $(I\times F,-d\tau^{2}+f^{2}g)$  with weakly convex fiber $(F,g)$ and satisfying Condition (A) is geodesically connected.
\end{lemma}

{\it Proof.}  Let $(\tau_{0},x_{0}), (\tau_{0}',x'_{0})\in I\times F$,  $L=d(x_{0},x'_{0})$, $L>0$ be.  We consider the following cases according to the values of $m_l , m_r$ in (\ref{e-m}):      

(i) Case $m_{l},m_{r}<\frac{1}{f^{2}(\tau_{0})}$. Then  $\int_{\tau_{0}}^{b_{\star}}f^{-2}(\frac{1}{f^{2}}-m_{r})^{-1/2}=\infty$      , $\int_{a_{*}}^{\tau_{0}}f^{-2}(\frac{1}{f^{2}}-m_{l})^{-1/2}=\infty$  and, thus, there exist $a_{*}<\tau^{-}<\tau_{0}<\tau^{+}<b_{*}$   such that $\int_{\tau_{0}}^{\tau^{+}}f^{-2}(\frac{1}{f^{2}}-m_{r})^{-1/2}=L$, $\int_{\tau^{-}}^{\tau_{0}}f^{-2}(\frac{1}{f^{2}}-m_{l})^{-1/2}=L$; so $(\tau_{0}, x_{0})$ can be joined with $(\tau_{\pm},x'_{0})$. By using   Lemma \ref{l.nuevo} we can connect $(\tau_{0},x_{0})$ with $[\tau^{-},\tau^{+}]\times \{x'_{0}\}$ taking $K\in [K^{-},K^{+}]$. Moreover, fixed $\epsilon >0$ such that $\tau^{+}+\epsilon<b$ (resp.  $\tau^{-}-\epsilon>a$) the limit of  $\int_{\tau_{0}}^{\tau^{+}+\epsilon}f^{-2}(\frac{1}{f^{2}}-D)^{-1/2}$ (resp.  $\int_{\tau^{-}-\epsilon}^{\tau_{0}}f^{-2}(\frac{1}{f^{2}}-D)^{-1/2}$) is greater than $L$ when $D\rightarrow m_{r}$ (resp.  $D\rightarrow m_{l}$) and the limit is $0$ when $D\rightarrow -\infty$; so, $(\tau_0,x_0)$ can be connected with $(\tau^+ + \epsilon, x_0')$ and 
$(\tau^- - \epsilon, x_0')$.                                                                                   Therefore, we can also connect $(\tau_{0},x_{0})$  with $[\tau^{+},b)\times \{x'_{0}\}$  and  $(a,\tau^{-}]\times \{x'_{0}\}$ taking $K\in[K^{+},\infty)$ and $K\in (-\infty,K^{-}]$, respectively. In particular, $(\tau_{0},x_{0})$, $(\tau'_{0},x'_{0})$ can be joined. 

(ii) Case $m_{l}=m_{r}=\frac{1}{f^{2}(\tau_{0})}$. Assume, say, $\tau_{0}<\tau'_{0}$;  then  $\int_{\tau_{0}}^{\tau'_{0}}f^{-2}(\frac{1}{f^{2}}-D)^{-1/2}$ goes to $0$ if $D\rightarrow -\infty$ and to $\infty$ if $D\nearrow \frac{1}{f^{2}(\tau_{0})}$. Therefore, there exist $D^{*}<\frac{1}{f^{2}(\tau_{0})}$ such that  $\int_{\tau_{0}}^{\tau'_{0}}f^{-2}(\frac{1}{f^{2}}-D^{*})^{-1/2}=L$ and the proof is over.  

(iii) Case  $m_{l}=\frac{1}{f^{2}(\tau_{0})}$ and $m_{r}<\frac{1}{f^{2}(\tau_{0})}$ (the remaining case is analogous). If, for certain $\delta>0$,  $\int_{a}^{\tau_{0}}f^{-2}(\frac{1}{f^{2}}-m_{l}-\delta)^{-1/2}>L$  then we can follow an argument as in (i). Otherwise, let $\tau^{+}$ be such that  $\int_{\tau_{0}}^{\tau^{+}}f^{-2}(\frac{1}{f^{2}}-m_{r})^{-1/2}=L$. Fixed $\epsilon >0$, the limit of  $\int_{\tau_{0}}^{\tau^{+}+\epsilon}f^{-2}(\frac{1}{f^{2}}-D)^{-1/2}$ is $0$ when $D\rightarrow -\infty$ and it is greater than $L$ when  $D\rightarrow m_{r}$; thus, we can connect  $(\tau_{0},x_{0})$ with $(\tau^+ + \epsilon, x_0')$ and, therefore, with  $[\tau^{+},b)\times \{x'_{0}\}$, by means of geodesics with $K\in[K^{+},\infty)$. Finally, from Lemma \ref{l1.75}, we can also connect $(\tau_{0},x_{0})$ with  $(a,\tau^{+}]\times \{x'_{0}\}$ taking $K\in(0,K^{+}]$. $\Box$

 \begin{lemma}\label{l3} A GRW spacetime  $(I\times F,-d\tau^{2}+f^{2}g)$  with weakly convex fiber $(F,g)$ and satisfying Condition (B) is geodesically connected.
\end{lemma}

{\it Proof.} Let $(\tau_{0},x_{0}), (\tau_{0}',x'_{0})\in I\times F$,  $L=d(x_{0},x'_{0})$, $L>0$ be. Firstly, suppose  the case 
\begin{equation}
\label{casofeo}
\int_{\tau_{0}}^{b_{\star}}f^{-2}(\frac{1}{f^{2}}-m_{r})^{-1/2}\leq L\quad \int_{a_{\star}}^{\tau_{0}}f^{-2}(\frac{1}{f^{2}}-m_{l})^{-1/2}\leq L
\end{equation}
$m_{r}=m_{b}$, $m_{l}=m_{a}$ and $2d_{b}\geq L$, $2d_{a}\geq L$. 

Note that 
\begin{equation} \label{maximo}
\frac{1}{f^2(\tau_0)} \geq Max\{ m_a , m_b \},
\end{equation}
and we consider first that this inequality is strict. 
Then, fixed $\delta>0$ such that $a+\delta<\tau_{0}, \tau'_{0}$ and $\tau_{0},\tau'_{0}<b-\delta$, there exist $0<K_{\delta}^{r}<\frac{1}{f^{2}(\tau_{0})}-m_{b}$ and $m_{a}-\frac{1}{f^{2}(\tau_{0})}<K_{\delta}^{l}<0$ such that $\tau(K_{\delta}^{r})>b-\delta$ and  $\tau(K_{\delta}^{l})<a+\delta$; recall that, otherwise, say
$$\begin{array}{c}
2 \int_{b-\delta}^{b_{*}}f^{-2}(\frac{1}{f^{2}}-D)^{-1/2}<L\leq 2d_{b}= \\ 2 lim\; sup_{\hat D\rightarrow m_{b}}(\int_{b-\delta}^{b_{*}}f^{-2}(\frac{1}{f^{2}}- \hat D)^{-1/2})-2\int_{b-\delta}^{b}f^{-2}(\frac{1}{f^{2}}-m_{b})^{-1/2}\end{array}$$
for all $D>m_{b}$ (with $b_*(D) > b-\delta)$, which is a contradiction because $\int_{b-\delta}^{b}f^{-2}(\frac{1}{f^{2}}-m_{b})^{-1/2}>0$.

So, the geodesics corresponding to $K_{\delta}^{r}$ and $K_{\delta}^{l}$ join $(\tau_{0},x_{0})$ with $(\tau^{r},x'_{0})$ and $(\tau^{l},x'_{0})$, where $\tau^{r} = \tau(K_{\delta}^{r})$,  $\tau^{l} = \tau(K_{\delta}^{l})$. From  Lemma \ref{l.nuevo}  we can connect $(\tau_{0},x_{0})$  with  $[\tau^{l},\tau^{r}]\times\{x'_{0}\}$ taking $K\in [K_{\delta}^{l},K_{\delta}^{r}]$ and, thus, the connectedness of $(\tau_{0},x_{0})$ with $(\tau'_{0},x'_{0})$ is obtained.

If (\ref{maximo}) holds with equality, then, because of (\ref{casofeo}) we have $m_a \neq m_b$ (say, $m_a > m_b$), and $K=0$ does not belong to the domain ${\cal D}$ of $\tau(K)$. Reasoning as above $K^+ \in {\cal D}, K^+>0$ is found, and the result follows from Lemma \ref{l1.75}.  

Finally, the reimaning cases (where not necessarily both inequalities (\ref{casofeo}) hold) are combinations of this one and the cases in Lemma \ref{l2}. $\Box$

\vspace{3mm}
\noindent Now, we are ready to prove our main result on connectedness.
The proof of (ii) in Theorem \ref{t0} is the consequence of the following two Propositions.

\begin{proposition}\label{p1} Let $(I\times F,-d\tau^{2}+f^{2}g)$ be a GRW spacetime with weakly convex fiber $(F,g)$ and satisfying Condition (C). Then it is geodesically connected.
\end{proposition}

{\it Proof.} Let $(\tau_{0},x_{0}), (\tau_{0}',x'_{0})\in I\times F$,  $L=d(x_{0},x'_{0})$, $L>0$ be. Suppose 
$$\int_{\tau_{0}}^{b_{\star}}f^{-2}(\frac{1}{f^{2}}-m_{r})^{-1/2}\leq L\quad \int_{a_{\star}}^{\tau_{0}}f^{-2}(\frac{1}{f^{2}}-m_{l})^{-1/2}\leq L$$ $m_{a}=m_{l}<m_{r}=m_{b}$, $2d_{b}<L \leq diam(F)$, $2d_{a}\geq L$   and $d_{b}\geq i_{b}$ (from Lemma \ref{l3} this is the only relevant case to study). As $d_{b}\geq i_{b}$  there exist $D_{1}^{r}\leq m_{b}$  such that  $$\int_{a_{*}}^{\tau_{0}}f^{-2}(\frac{1}{f^{2}}-D_{1}^{r})^{-1/2}+\int_{a_{*}}^{b}f^{-2}(\frac{1}{f^{2}}-D_{1}^{r})^{-1/2}<2d_{b}<L.$$ 
On the other hand, as $2d_{a}\geq L$, for $D_{2}^{r}<D_{1}^{r}$ near enough to $m_{l}$ we have $$\int_{a_{*}}^{\tau_{0}}f^{-2}(\frac{1}{f^{2}}-D_{2}^{r})^{-1/2}+\int_{a_{*}}^{b}f^{-2}(\frac{1}{f^{2}}-D_{2}^{r})^{-1/2}>L.$$ 
Therefore,  the domain ${\cal D}$ of $\tau(K)$ contains $K_{2}^{r}=D_{2}^{r}-\frac{1}{f^{2}(\tau_{0})}$ but not  $K_{1}^{r}=D_{1}^{r}-\frac{1}{f^{2}(\tau_{0})}$. From Lemma \ref{l1.75}, $(\tau_0,x_0)$ can be connected with $[\tau(K_2^r),b) \times \{ x'_0 \}$. Choose $K^r= D^{r}-\frac{1}{f^{2}(\tau_{0})}
\in [K_{2}^{r}, K_{1}^{r})$ such that $\tau(K^{r})  
>b-\delta$, for $\delta$ small. 
As $2d_{a}\geq L$, there exist $D^l$, $m_{a}<D^{l}<D^{r}$ such that $\tau(K^{l})<a+\delta$ ($K^{l}=D^{l}-\frac{1}{f^{2}(\tau_{0})}$). Thus, the result follows from Lemma \ref{l.nuevo}. 
$\Box$

\begin{proposition}\label{p2} Let $(I\times F,-d\tau^{2}+f^{2}g)$ be a GRW spacetime with weakly convex fiber $(F,g)$ and satisfying Condition (R). Then it is geodesically connected.
\end{proposition}

{\it Proof.} We will use sistematically that if $D$ is close enough to $m$ and $D>m$ then $K^+ = 1/f^2(\tau_0) -D (>0)$ and $K^- = D-1/f^2(\tau_0) (<0) $ satisfy $[K^-,K^+] \subset {\cal D}$; thus, Lemma \ref{l1.5} can be claimed. Let $(\tau_{0},x_{0}), (\tau_{0}',x'_{0})\in I\times F$,  $L=d(x_{0},x'_{0})$, $L>0$ be, and consider the following two cases:

(i) Suppose  $r_{i}^{n}(\tau_{0})\leq L \leq r_{s}^{n}(\tau_{0})$,  $l_{i}^{n'}(\tau_{0})\leq L \leq l_{s}^{n'}(\tau_{0})$ for certain $n,n'\geq 0$. Fix $\epsilon >0$ such that $a+\epsilon <\tau'_{0}<b-\epsilon$. Then for some $D_{i}^{r}, D_{s}^{r}$ close to $m$, chosen such that $m< D_{i}^{r} < D_{s}^{r}$, we have 

\begin{equation}
\begin{array}{l}
_{(-)^{n}[n-1]}\int_{\tau_{0}}^{a_{*}(D_{i}^{r})}f^{-2}(\frac{1}{f^{2}}-D_{i}^{r})^{-1/2}+\int_{a_{*}(D_{i}^{r})}^{b-\epsilon}f^{-2}(\frac{1}{f^{2}}-D_{i}^{r})^{-1/2}<L \\ _{(-)^{n}[n]}\int_{\tau_{0}}^{b_{*}(D_{s}^{r})}f^{-2}(\frac{1}{f^{2}}-D_{s}^{r})^{-1/2}+\int_{b-\epsilon}^{b_{*}(D_{s}^{r})}f^{-2}(\frac{1}{f^{2}}-D_{s}^{r})^{-1/2}>L
\end{array}
\end{equation}
if $n\geq 1$, or
\begin{equation}
\begin{array}{l}
\int_{\tau_{0}}^{b-\epsilon}f^{-2}(\frac{1}{f^{2}}-D_{i}^{r})^{-1/2}<L \\ \int_{\tau_{0}}^{b_{*}(D_{s}^{r})}f^{-2}(\frac{1}{f^{2}}-D_{s}^{r})^{-1/2}+\int_{b-\epsilon}^{b_{*}(D_{s}^{r})}f^{-2}(\frac{1}{f^{2}}-D_{s}^{r})^{-1/2}>L
\end{array}
\end{equation}
if $n=0$. Reasoning similarly to the left, we obtain analogous $D_{i}^{l}, D_{s}^{l}$, with $m<D_{i}^{l}<D_{s}^{l}$. From Lemma \ref{l1} there exist $D^r , D^l$, with $D_{i}^{r}<D^{r}<D_{s}^{r}$, $D_{i}^{l}<D^{l}<D_{s}^{l}$ such that $\tau(K^{r})>b-\epsilon$, $\tau(K^{l})<a+\epsilon$, where $K^{r}=(-1)^{n} (\frac{1}{f^{2}(\tau_{0})}-D^{r})$, $K^{l}=(-1)^{n'-1}(\frac{1}{f^{2}(\tau_{0})}-D^{l})$. Therefore, as $a+\epsilon<\tau'_{0}<b-\epsilon$, the connectedness of  $(\tau_{0},x_{0})$ with $(\tau'_{0},x'_{0})$ is a consequence of  Lemma \ref{l1.5}.

(ii) Suppose now $L<r_{i}^{0}(\tau_{0})(<r_{i}^{n}(\tau_{0}))$ and $L<l_{i}^{0}(\tau_{0})(<l_{i}^{n}(\tau_{0}))$. As we saw below Definition \ref{d2.5},  $r_{i}^{0}(\tau_{0})= \int_{\tau_{0}}^{b}f^{-2}(\frac{1}{f^{2}}-m)^{-1/2}$ (analogoussly for $l_i^0$), thus,  there exist  $\epsilon >0$, $a+\epsilon<\tau'_{0}<b-\epsilon$ such that 
$$\int_{\tau_{0}}^{b-\epsilon}f^{-2}(\frac{1}{f^{2}}-m)^{-1/2}>L, \quad 
\int^{\tau_{0}}_{a+\epsilon}f^{-2}(\frac{1}{f^{2}}-m)^{-1/2}>L.$$ 
But the limit of $\int_{\tau_{0}}^{b-\epsilon}f^{-2}(\frac{1}{f^{2}}-D)^{-1/2}$ when $D\rightarrow -\infty$ is $0$, thus we obtain $D^{r}<m$ such that $\int_{\tau_{0}}^{b-\epsilon}f^{-2}(\frac{1}{f^{2}}-D^{r})^{-1/2}=L$. So, taking $K^{r}=\frac{1}{f^{2}(\tau_{0})}-D^{r} (>0)$, we obtain $\tau(K^{r})=b-\epsilon$. Analogously, there exist $K^{l}<0$ such that $\tau(K^{l})=a+\epsilon$. Therefore,  we obtain the connectedness of $(\tau_{0},x_{0})$ with $(\tau'_{0},x'_{0})$ from Lemma \ref{l1.5} again. 

The remaining cases are combinations of the previous ones. $\Box$

\vspace{3mm}

\noindent {\it Proof of Theorem \ref{t1}.} For (i) assume that  $z_{0}=(\tau_{0},x_{0})$, $z'_{0}=(\tau'_{0},x'_{0})$ are causally related and $\tau_0 < \tau'_0$. From Theorem \ref{t0} there exist a non-spacelike geodesic $\gamma:{\cal J}\rightarrow I\times F$, $\gamma(t)=(\tau(t),\gamma_{F}(t))$ joining them. As $(F,g)$ is strongly convex, necessarily

\begin{equation}
\int_{\tau_{0}}^{\tau'_{0}}f^{-2}(\frac{1}{f^{2}}-D_{0})^{-1/2}=d(x_{0},x'_{0})
\end{equation}
being $D_{0}=g(\frac{d\gamma}{dt},\frac{d\gamma}{dt})\leq 0$. But the integral $\int_{\tau_{0}}^{\tau'_{0}}f^{-2}(\frac{1}{f^{2}}-D)^{-1/2}$ is strictly increasing with $D$, for $D\leq 0$; thus, $\gamma$ is the only causal geodesic joining $z_{0}$ and $z'_{0}$. Moreover, when $D>0$ the integral (possibly under Convention \ref{cv1}) is bigger than when $D=0$; so, no spacelike geodesic joins $z_{0}$ and $z'_{0}$. 

In order to prove (ii) assume that neither Condition (C) nor Condition (R) hold and consider the following cases. In the first three ones we will assume that Condition (R) is not appliable, and Condition (C) does not hold at $b$ (at $a$ would be analogous). Recall that, from Lemma \ref{l00}, $1/f^2$ is decreasing at $b$; in the first case $b$ is a non-unique absolute minimum; in the second, $b$ is the unique absolute minimum, which is simpler; in the third, $b$ is not an absolute minimum, which oblies to use properly the definition of $i_b$. In the fourth case, Condition (R) is appliable, but it does not hold (neither does Condition (C), see Remark \ref{r1}).

(i) Assume that $b$ is a relative minimum of $1/f^2$ and $m=m_{b}$ is reached at a point $\tau_{m}\in (a,b)$. As $2d_{b}<diam(F)$, choose $L>0$  such that $2d_{b}<L<diam(F)$. From this choice,  there exist $\tau_{0}^{r} 
> \tau_{m}$, close to $b$ such that 
\begin{equation} \label{nogeo1}
2\int_{\tau_{0}^{r}}^{b_{*}(D)}f^{-2}(\frac{1}{f^{2}}-D)^{-1/2}<L, \quad  \forall D\in (m,\frac{1}{f^{2}(\tau_{0}^{r})}). 
\end{equation}
As   $\tau_{m}$ is a minimum,  $\frac{d}{d\tau}\frac{1}{f^{2}}\mid_{\tau_{m}}=0$. Thus, there exist $\tau_{0}^{l}$ near enough to $b$ such that 
\begin{equation} \label{nogeo2}
2\int_{a_{\star}(D)}^{\tau_{0}^{l}}f^{-2}(\frac{1}{f^{2}}-D)^{-1/2}>L \quad \forall D\in (m,\frac{1}{f^{2}(\tau_{0}^{l})}). 
\end{equation}
Now, taking any $\tau_{0}  > Max\{\tau_{0}^{r},\tau_{0}^{l}\}$, $\tau'_0 > \tau_0$ and $x_0, x'_0$ with $d(x_0, x'_0)= L$, it is clear that (\ref{nogeo1}) and (\ref{nogeo2}) forbid to connect $(\tau_{0},x_{0})$, $(\tau'_{0},x'_{0})$  by means of a geodesic.

(ii) Assume that $b$ is a relative minimum, $m=m_{b}<m_{a}$ and $\frac{1}{f^{2}(\tau)}>m$ for all $\tau\in (a,b)$. Then, necessarily, $2d_{b}<diam(F)$. Choose again $\tau_{0}^{r}$ such that (\ref{nogeo1}) hold. Recall that we can impose now, aditionally, that $\tau_{0}^{r}$ is the strict minimum of $1/f^2$ on $(a,\tau_{0}^{r}]$. So, clearly, $(\tau_{0}^{r},x_{0})$ and $(\tau'_{0},x'_{0})$ cannot be joined by a geodesic, if $\tau_{0}^{r} < \tau'_0$ and $d(x_0,x'_0)=L$.

(iii) Assume that $b$ is a relative minimum and $m<m_{b}$. As $d_{b}<i_{b}$  there exist $\tau_{0}^{l}$ such that 
\begin{equation} \label{nogeo3}
2\int_{a_{*}(D)}^{\tau_{0}^{l}}f^{-2}(\frac{1}{f^{2}}-D)^{-1/2}\geq 2d_{b}+ 2\epsilon, \quad \forall D\in (m,m_{b}]
\end{equation}
for some $\epsilon>0$ such that $2d_{b}+ 2\epsilon<diam(F)$. From the continuity stated in Lemma \ref{l1}, there exist $\delta >0$ such that 
inequality (\ref{nogeo3}) holds if the right member is replaced by $L= 2d_{b}+\epsilon $,  for all $D\in (m,m_{b}+\delta ]$. 

Now, as in case (i) we can take $\tau_{0} (=\tau_0^r) >\tau_{0}^{l}$, with $\frac{1}{f^{2}(\tau_{0})}<m_{b}+\delta$ and such that (\ref{nogeo1}) holds for all $D$. Thus, for any $\tau'_{0}>\tau_{0}$ we cannot connect $(\tau_{0},x_{0})$, $(\tau'_{0},x'_{0})$ by means of a geodesic, if  
 $d(x_{0},x'_{0})=L$.

(iv) Assume that $\frac{1}{f^{2}(\tau)}>m$ for $\tau \in (a,b)$ and $m_{a}=m_{b}=m$. Suppose that Condition (R) is not fulfilled by, say, the $r'$s, that is, $r_{s}^{n}(\tau_{0})<r_{i}^{n+1}(\tau_{0})$ with $r_{s}^{n}(\tau_{0})<diam(F)$ for certain $n\geq 0$ and $\tau_{0}\in I$ (see the comments below Definition \ref{d2.5}). Fix $L=d(x_{0},x'_{0}) (\leq diam(F))$ with  $r_{s}^{n}(\tau_{0})<L<r_{i}^{n+1}(\tau_{0})$. These inequalities imply for $n\geq 1$ that 
there exist $\epsilon >0$ such that 

\begin{equation}
\begin{array}{l}lim\; inf_{D\searrow m}\{_{(-)^{k}[k-1]}\int^{a_{*}}_{\tau_{0}}f^{-2}(\frac{1}{f^{2}}-D)^{-1/2}+\int_{a_{*}}^{b-\epsilon}f^{-2}(\frac{1}{f^{2}}-D)^{-1/2}\}>L \\ lim\; sup_{D\searrow m}\{_{(-)^{k'}[k']}\int_{\tau_{0}}^{b_{*}}f^{-2}(\frac{1}{f^{2}}-D)^{-1/2}+\int_{b-\epsilon}^{b_{*}}f^{-2}(\frac{1}{f^{2}}-D)^{-1/2}\}<L \end{array}
\end{equation}  
for $k=n+1$, $k'=n$ and, thus, for all $k\geq n+1$ and $k'\leq n$. But this implies that, for some $\delta >0$ with $b^{*}(D=m+\delta)>b-\epsilon$, if $m<D<m+\delta$ then

\begin{equation}
\begin{array}{l}_{(-)^{k}[k-1]}\int^{a_{*}}_{\tau_{0}}f^{-2}(\frac{1}{f^{2}}-D)^{-1/2}+\int_{a_{*}}^{b-\epsilon}f^{-2}(\frac{1}{f^{2}}-D)^{-1/2}>L \\ _{(-)^{k'}[k']}\int_{\tau_{0}}^{b_{*}}f^{-2}(\frac{1}{f^{2}}-D)^{-1/2}+\int_{b-\epsilon}^{b_{*}}f^{-2}(\frac{1}{f^{2}}-D)^{-1/2}<L \end{array}
\end{equation}  
for $k\geq n+1$, $k'\leq n$ (there are analogous inequalities when $n=0$). Therefore,  $(\tau_{0},x_{0})$ cannot be geodesically connected with $(\tau'_{0},x'_{0})$ if $\tau'_{0}>b-\epsilon$. $\Box$

\vspace{3mm}

{\it Discussion.} Next, we will see what happens if we assume just weak convexity in Theorem \ref{t1} and Condition (R) is appliable (a similar study could be done if Condition (C) is appliable instead). As a consequence, we will give a proof of the (well-known) non-geodesic connectedness of de Sitter spacetime. It should be noticed that previous proofs use the high degree of symmetry of this spacetime \cite{CM}, \cite{Sc}. In our proof we will see what is the exact role of this symmetry.

Fix $z_0 = (\tau_0 , x_0) \in I\times F, x'_0 \in F, x_0 \neq x'_0$ and $\epsilon >0$. Put 
$r_{i, \epsilon}^{0}(\tau_{0}), r_{s, \epsilon}^{n}(\tau_{0})$ etc. equal to the quantities in Definition \ref{d2.5} but {\it without} taking the limit $\epsilon \rightarrow 0$ (the extension of this new definition when $a= -\infty $ or $b= \infty$ is obvious, see de Sitter spacetime below). Now, consider: 
$$A_{\epsilon}=  \cup_{n\geq 0}[r_{i, \epsilon}^{n}(\tau_{0}),r_{s, \epsilon}^{n}(\tau_{0})] \cup [0,r_{i, \epsilon}^{0}(\tau_{0})], \; \quad  \;      
B_{\epsilon}=  \cup_{n\geq 0}[l_{i,\epsilon}^{n}(\tau_{0}),l_{s,\epsilon}^{n}(\tau_{0})] \cup 
[0,l_{i,\epsilon}^{0}(\tau_{0})],$$
 and also
$${\cal L} = \{ \mbox{length}(\hat \gamma_F) \mid \hat \gamma_F \; \mbox{is a F-geodesic which joins} \; x_0 \; \mbox{and} \; x'_0 \} \subset (0,\infty).$$
From the proof of Theorem \ref{t1} $z_0$ can be joined with $[a+\epsilon_0,b-\epsilon_0] \times \{ x_0'\}$ if 
$$
{\cal L} \cap A_{\epsilon} \cap B_{\epsilon} \neq \emptyset,
$$
for some $\epsilon < \epsilon_0$. Moreover, it is also clear that $z_0$ {\em cannot} be joined with the points in $(a, a+\epsilon_0) \times \{x_0'\} \; \cup 
\; (b-\epsilon_0,b) \times \{x_0'\} $ if 
\begin{equation}\label{noconge}
{\cal L} \cap A_{\epsilon_0} = \emptyset \quad \mbox{and} \quad 
{\cal L} \cap  B_{\epsilon_0} = \emptyset.
\end{equation}

For de Sitter spacetime, $I=\R{}$, $f=\cosh$ and the fiber is the usual sphere of radius 1. Recall that when the interval $I$ is not bounded, we must replace  $ b-\epsilon$ (if $b=\infty$) and  $-(a+\epsilon)$ (if $a=-\infty$) by $M>0$, and the limit $\epsilon \rightarrow 0$ must be replaced by $M\rightarrow \infty$. 
Take $z_0= (0, x_0)$;  by Definition \ref{d2.5} ($M \rightarrow \infty$) we have:
\begin{equation} \label{m8}
 r_{i}^{n}(0) = r_{s}^{n}(0) = \frac{\pi}{2} + n\pi 
= 
 l_{i}^{n}(0)=l_{s}^{n}(0).
\end{equation}
For $M=0$, the new definitions $r_{i, \epsilon}^{0}(0), r_{s, \epsilon}^{n}(0) \; (\epsilon\equiv \infty)$ reads:
\begin{equation}\label{m0}
\begin{array}{l}
 r_{i, \epsilon}^{n}(0) = n\pi = l_{i,\epsilon}^{n}(0) \\
 r_{s, \epsilon}^{n}(0) =  (n+1) \pi =
l_{s,\epsilon}^{n}(0). 
\end{array}
\end{equation}
Now, choose $x'_0$ as the antipodal point of $x_0$, that is:
$${\cal L} = \{(2n+1) \pi \mid n=0,1,2\dots \}.$$
From the two limit cases (\ref{m8}), (\ref{m0}), it is clear that condition (\ref{noconge}) is fulfilled  for any $M>0$. So $z_0$ {\it cannot} be joined by means of a geodesic with $(-\infty, 0) \times \{x_0'\} \; \cup \;
(0,\infty) \times \{x_0'\}$.

Summing up, for de Sitter spacetime the ``symmetries'' of its warping function  are essential in order to have enough ``holes'' in $A_{\epsilon_0}$ and $B_{\epsilon_0}$, where all the elements of ${\cal L}$ lie. But the only relevant symmetry of the fiber is that there are two points $x_0, x'_0$ such that the lengths of the geodesics which joins them has a constant gap. In our case, this gap ($2\pi$) and the symmetries of $f$ fits well when $d(x_0, x'_0)=\pi$.

\vspace{3mm}

\noindent {\it Proof of Theorem \ref{t2}.} For the first assertion consider a F-geodesic $\hat \gamma_{F}(r)$, $r\in [0,L]$ with $L=long \hat \gamma_{F}$, $\hat \gamma_{F}(0)=x_{0}$ and $\hat \gamma_{F}(L)=x'_{0}$. From our hypotheses, if $\frac{1}{f^{2}}$ reaches a relative minimum at $b$ (resp. $a$) and $b^*(m_r)=b$ (resp. $a^*(m_l)=a$) then  

\begin{equation} \label{fff}
\begin{array}{l}
\mbox{either} \quad \int^{b^*(m_r)}_{\tau_{0}}f^{-2}(\frac{1}{f^{2}}-m_{r})^{-1/2}>L \quad \mbox{or} \quad 2d_{b}>L \\

(\mbox{resp.  either} \quad \int_{a^*(m_l)}^{\tau_{0}}f^{-2}(\frac{1}{f^{2}}-m_{l})^{-1/2}>L \quad \mbox{or} \quad 2d_{a}>L). 
\end{array}
\end{equation}
As we checked in Lemma \ref{l2}  and Lemma \ref{l3},  inequalities (\ref{fff}) allow us to obtain a geodesic joining $z_{0}$ and $z'_{0}$ with component on the fiber a reparameterization of $\hat \gamma_{F}(r)$ (recall that in these lemmas $\hat \gamma_{F}(r)$ was always taken as a minimizing $F-$geodesic, but the minimizing property was used just to ensure that (\ref{fff}) hold). It is straightforward to check that these inequalities also hold if $b^*(m_r)<b$ or $a< a^*(m_l)$, because the corresponding integral is then infinite.

If $(F,g)$ is complete  and $F$ is not contractible then, fixed $x_{0}, x'_{0} \in  F$,  there exist a sequence of  geodesics $\hat \gamma^m_{F}(r)$ 
joining $x_{0}$ and $x'_{0}$, with diverging lengths $L_m$ (see for example \cite[Th.2.11.9]{Ma}). 
Let $\gamma^m(t)= (\tau^m(t), \gamma^m_{F}(t))$ be the geodesic connecting $z_0, z'_0$ constructed from $\hat \gamma^{m}_{F}(r)$, 
and assume  $\tau_0 \leq  \tau'_0$.  If $\gamma^m(t)$ is causal, then necessarily (\ref{e-8}), (\ref{h}) hold with $L=L_m$. But in this case $D \leq 0$ and, thus, $L_m \leq \int_{\tau_{0}}^{\tau'_{0}}\frac{1}{f} (<\infty)$.
As the sequence $\{L_m\}$ is diverging, all the geodesics but a  finite number are spacelike.

The last assertion is also a direct consequence of the  fact that the lengths of the $F-$pregeodesics corresponding to causal geodesics are  bounded by $\int_{\tau_{0}}^{\tau'_{0}}\frac{1}{f}$, and Lemma \ref{l100}. $\Box$ 

\begin{lemma}\label{l100} If $(M,g)$ is a complete Riemannian manifold and $p,q\in M$ are no conjugate then for all $L>0$ there exist at most finitely many geodesics with length smaller than $L$ connecting $p$ and $q$.
\end{lemma}

{Proof.} Otherwise, from the compactness of $\{v\in T_{p}M: \mid v \mid \leq L\}$, we would obtain a sequence $\{v_{n}\}_{n\in \N{}}$, $v_{n}\rightarrow v_{0}$, $v_{n},v_{0}\in T_{p}M$ such that $exp_{p}(v_{n})=exp_{p}(v_{0})=q$ for all $n$. Then, $v_{0}$ would be a singular point of $exp_{p}$ and, thus, $p$ and $q$ would be conjugate for the geodesic $\gamma(t)=exp_{p}(t\cdot v_{0})\quad t\in \R{}$,  which is a contradiction. $\Box$

{\it Remark.} In the proof of Theorem \ref{t2} we have used that, for a complete Riemannian manifold which is non--contractible in itself, infinitely many geodesics joining $p$ and $q$ exist, and there is a sequence of them with diverging lengths. So, in this case, Lemma \ref{l100} says that if $p$, $q$ are not conjugate then any sequence of geodesics joining them have diverging lengths. In particular, the number of geodesics joining two non-conjugate points of a complete Riemannian manifold must be enumerable.

\section{Conjugate points and Morse-type inequalities.} \label{s5}
\label{s1}

In order to prove results on conjugate points, it seems more natural to consider all the geodesics obtained by varying a fixed one with the same speed $D$. So, we will drop previous normalization $c=1$ for geodesics non-tangent to the base. The only modification in previous formulae which we will have to bear in mind is that, now, (\ref{e-5}) reads
\begin{equation}\label{e-5bis}
\frac{dt}{dr}=
\frac{1}{\sqrt{c}} \cdot f^{2}\circ\tau\circ t
\end{equation}
so, the definition of $h$ in (\ref{h}) must be changed to 
\be
\label{hc}
h^{\epsilon}=\epsilon \sqrt{c} \cdot f^{-2}(-D+\frac{c}{f^{2}})^{-1/2}.
\ee

\begin{theorem} 
\label{t-3}
Let  $z_0=(\tau_{0},x_{0})$, $z'_0=(\tau'_{0},x'_{0})$ be two points of the GRW   spacetime $(I\times F,-d\tau^{2}+f^{2}g)$  with $n$--dimensional fiber $(F,g)$. Assume that $\gamma(t)=(\tau(t),\gamma_{F}(t))$ is a geodesic which joins them, being $\gamma_{F}(t)$ the reparameterization of a non-constant $F$-geodesic $\hat \gamma_{F}$, and that $z_0, z_0'$ are conjugate along $\gamma$ with multiplicity $m\in \{0, 1 \dots n\}$ ($m=0$ means no conjugate). 

(i) Then $x_{0}$, $x'_{0}$ are conjugate points of multiplicity $m'\in \{m,m-1\}$  along $\hat{\gamma}_{F}$ (at the corresponding points of the domain). In particular, if $z_0, z'_0$ are non-conjugate then so are $x_0$ and $x_0'$.

(ii) If $\gamma $ is a causal geodesic (or any geodesic without zeroes in $d\tau/dt$) then $m'=m$.
\end{theorem}
 
\noindent {\it Remark.} (1) The following direct computation shows that, even in the excluded case $\hat \gamma_{F} \equiv x_0=x_0'$ ($\hat \gamma_{F}$ is constant), the points  $z_0 , z'_0$ are not conjugate. Thus this case can be  included in Theorem \ref{t-3} with the convention ``a constant geodesic 
$\hat \gamma_F$ has no conjugate points''. Assume $\tau_0 < \tau'_0$ and consider the geodesic $\gamma(t)= (t,x_0), \quad t\in[\tau_0,\tau'_0]$. Let $E_i(t), i\in\{1,\dots n\}$ be orthonormal parallel fields along $\gamma$ which span the orthogonal to $\gamma'$. A vector field $J(t)= \sum_i a_i(t)E_i(t)$ along $\gamma$ is a Jacobi field if and only if each function $a_i(t)$ is a solution of the Sturm differential equation:
\begin{equation} \label{je}
a''(t) - \frac{f''(t)}{f(t)} a(t) = 0, \quad \quad  t\in[\tau_0,\tau'_0].
\end{equation}
But, clearly, $f(t)$ is also a strictly positive solution of (\ref{je}). Thus, if $a(\tau_0)=0$ and $a'(\tau_0)\neq 0$ then $a(\tau)$ cannot vanish on $(\tau_0,\tau'_0]$, as required.

(2) Moreover, for any $\tau > \tau_0$, replace (\ref{je}) by the spectral equation (see \cite{BGM}):
\begin{equation} \label{jo}
a''(t) - \frac{f''(t)}{f(t)} a(t) + \lambda_{\tau} a(t) = 0,  
\end{equation}
$\lambda_{\tau} \in \R{}$, with boundary conditions $a(\tau_0)=a(\tau)=0$. A simple Sturm argument shows that if $\tau < \bar \tau$ then $\lambda_{\tau} > \lambda_{\bar \tau}$; that is, the spectral flow $\lambda(\tau)\equiv\lambda_{\tau}$ is decreasing. This also holds for the static bidimensional case (see the next section), and  should be compared with  \cite{BGM}. At any case, the main result of \cite{BGM} can be reobtained,  as we will see in the next section; independently, it is also reobtained in \cite{GMPT}, in the general setting of geodesics admitting a timelike Jacobi field.

{\it Proof of Theorem \ref{t-3}.}  Step 1. For any geodesic $\gamma$, $m'\geq m-1$.

\noindent Consider $v_{0},v_{1},\ldots,v_{m}\in T_{z_{0}}(I\times F)$ such that $V=Span\{v_{1},\ldots, v_{m}\}$  where $V=ker((dexp_{z_{0}})_{v_{0}})$ and $exp_{z_{0}}(v_{0})=z'_{0}$. From semi-Riemannian Gauss Lemma \cite[5.1]{O} $v_{0}$ and each $v_{i}$ are orthogonal, so, $\{v_{0},\ldots,v_{m}\}$ are linearly independent (recall that if $v_{0}$ is lightlike then as $v_{0}$ and $v_{i}$ are not collinear then each $v_{i}$ is spacelike). Moreover, consider the usual projection on the fiber, $\pi_{F}$; as 
 $\gamma$ is not on the base, then $(d\;\pi_{F})_{z_{0}}v_{0}\neq 0$ and,  say, $\{(d\;\pi_{F})_{z_{0}}v_{0},(d\;\pi_{F})_{z_{0}}v_{1},\ldots ,
(d\;\pi_{F})_{z_{0}}v_{m-1}\}$ are linearly independent. So,  $(d\;\pi_{F})_{z_{0}}v_{0}$ is parallel to the initial velocity of $\hat \gamma_F$, and we have just to prove that there exist a direction of conjugacy of $\hat \gamma_F$ between $x_0, x'_0$ in each plane 
$W_i = Span\{(d\;\pi_{F})_{z_{0}}v_{0}, (d\;\pi_{F})_{z_{0}}v_{i}\} \subseteq T_{x_{0}}F$,  for $i = 1,\ldots, m-1$.

Defining  $\alpha_{i}(s)=v_{0}+sv_{i}$  we have  $\frac{d}{ds}\mid_{s=0}exp_{z_0}(\alpha_{i}(s))=0$ and, thus, 
\begin{equation} \label{nuc}
\frac{d}{ds}\mid_{s=0}\pi_{F}\circ exp_{z_0}(\alpha_{i}(s))=0. 
\end{equation}
There exist a non-constant continuous curve $\beta_{i}(s)\in W_i \quad i=1,\ldots ,m-1$ such that
\begin{equation}
\label{e1}
exp_{x_{0}}(\beta_{i}(s))=\pi_{F}\circ exp_{z_0}(\alpha_{i}(s)). 
\end{equation}
In fact, we take 
\begin{equation} \label{be}
\beta_i(s)= \mu_i(s) \frac{d\pi_F(\alpha_i(s))}{\mid d\pi_F(\alpha_i(s)) \mid} , 
\end{equation}
where $\mu_i(s)$ is the length of the pregeodesic 
$t\rightarrow \pi_{F}\circ exp_{z_0}(t\cdot \alpha_{i}(s))$ on $[0,1]$. 

Recall that  $(d\;\pi_{F})_{z_{0}}v_{0}$ is paralell to $\beta_i(0)\equiv \omega_0$, and
we had to prove that  $(dexp_{x_{0}})_{w_{0}}$ restricted to $W_i$ is singular. Otherwise, $\beta_{i}(s)$ would be smooth around 0 from (\ref{e1}). From (\ref{be}),  $0 \neq \beta'_{i}(0) \in W_i$, and from $(\ref{nuc})$ and $(\ref{e1})$, $\beta'_i(0) \in ker (dexp_{x_{0}})_{\omega_0}$, a contradiction.

Step 2. If $\gamma$ is causal then $m'\geq m$.

\noindent We will check that if $\gamma$ is not tangent to the base but it is causal (or any geodesic without zeroes in the derivative of the timelike component) then $\{v_{1},\ldots ,v_{m}\}$ are tangent to the fiber. So, $\{d\;\pi_{F})_{z_{0}}v_{0},(d\;\pi_{F})_{z_{0}}v_{1},\ldots ,
(d\;\pi_{F})_{z_{0}}v_{m}\}$ are linearly independent and the result follows as in previous step.

 From the hypotheses, 
\begin{equation}\label{e-400}\frac{d\tau}{dt}=\epsilon\sqrt{-D+\frac{c}{f^{2}\circ\tau}}\neq 0 \end{equation}
for all $t$ where $D=g^{f}(\frac{d\gamma}{dt},\frac{d\gamma}{dt})$ and $c=(f^{4}\circ \tau)\cdot g(\frac{d\gamma_{F}}{dt},\frac{d\gamma_{F}}{dt})$. Consider the usual projection on the base $\pi_{I}$, we will check that  $(d\pi_{I})_{z_{0}}(v_i)=0$. Let $\alpha_i(s)\in T_{z_{0}}(I\times F)$ be a curve such that $\alpha_i(0)=v_{0}$, $\frac{d}{ds}\mid_{s=0}\alpha_i(s)=v_i$, as above, and we also impose $g^{f}(\alpha_i(s),\alpha_i(s))=g^{f}(v_{0},v_{0})$ for all $s$. Put  $\gamma(s,t)=exp_{z_{0}}(t\cdot \alpha_i(s))\equiv(\tau_{s}(t),\gamma_{F\; s}(t))$ (thus $D(s)\equiv g^{f}(v_{0},v_{0})$). If $d\pi_{I}(v_i)=\frac{d}{ds}\mid_{s=0}\tau'_{s}(0)\neq 0$  then, as $D=-\tau'_{s}(0)^{2}+\frac{c(s)}{f^{2}(\tau_{0})}$ is constant, we obtain that $\frac{d}{ds}\mid_{s=0}c(s)\neq 0$. Now, including in (\ref{e-400}) the dependence on $s$ we have
$$\int_{\tau_{s}(0)}^{\tau_{s}(1)}\frac{d\tau}{\epsilon \sqrt{-D+\frac{c(s)}{f^{2}(\tau)}}}=1$$
and deriving with respect to $s$ we obtain $\frac{d}{ds}\mid_{s=0}\tau_{s}(1)\neq 0$. Therefore, $\frac{d}{ds}\mid_{s=0}\pi_{I}\circ exp_{z_{0}}(\alpha_i(s))\neq 0$ which contradicts that $v_i$ is a direction of conjugacy.

Step 3. $m\geq m'$.

\noindent Let $x_{0}$, $x'_{0}$ be conjugate points of multiplicity $m'$ along the F-geodesic $\hat{\gamma}_{F}$ and suppose $\gamma(0)=z_{0}$, $\gamma(1)=z'_{0}$. If $Span\{w_{1},\ldots ,w_{m'}\}=ker((d\; exp_{x_{0}})_{w_{0}})$ where $exp_{x_{0}}(w_{0})=x'_{0}$, consider a curve $\beta_{i}(s)$  in $T_{x_{0}}F$ such that $\beta_{i}(0)=w_{0}$, $\frac{d}{ds}\mid_{s=0}\beta_{i}(s)=w_{i}$ and $\mid \beta_{i}(s)\mid =\mid w_{0}\mid$ for all $s$, $i=1,\dots m'$. Define $\alpha_{i}(s)\in T_{z_{0}}(I\times F)$ such that $(d\pi_{I})_{z_{0}}(\alpha_{i}(s))= \frac{d\tau}{dt}(0)$ and 

\begin{equation}
\label{e4}
(d\pi_{F})_{z_{0}}(\alpha_{i}(s))=\frac{\sqrt{c}}{f^{2}(\tau_{0})}\cdot \frac{\beta_{i}(s)}{\mid w_{0}\mid }.
\end{equation}
For each $s$, the geodesic on the GRW spacetime $\gamma(s,t)=exp_{z_{0}}(t\cdot \alpha_{i}(s))\equiv(\tau_{s}(t),\gamma_{F\; s}(t))$, satisfy that $\gamma_{F \; s}(t)=exp_{x_{0}}(r_{s}(t)\cdot \beta_{i}(s))$ where $r_{s}(t)$ is an increasing function, because $\gamma_{F \; s}(t)$ is a pregeodesic on the fiber $F$. But, from (\ref{e-4}) and (\ref{e-5bis}), $r_{s}(t)$ is determined just by $c(s)\equiv c$ and $D(s)\equiv D$, so, $r_{s}(t)$ is independent of $s$, i.e. $r_{s}(t)\equiv r(t)$. Computing for $s=0$, it is clear that $r^{-1}(1)=1$ thus, necessarily $\pi_{F}\circ exp_{z_0}(\alpha_{i}(s))=exp_{x_{0}}(\beta_{i}(s))$ for all $s$. As $w_{i}\in ker((d\; exp_{x_{0}})_{w_{0}})$, we have: 
\begin{equation}\label{zz}
\frac{d}{ds}\mid_{s=0}\pi_{F}\circ exp_{z_0}(\alpha_{i}(s))=0.
\end{equation} 
On the other hand, from the relation between the parameters $\tau$ and $r$ for $\gamma(s,t)$ given by (\ref{hc})  we have 

\begin{equation}
\label{e5}
\int_{\tau_{0}}^{\tau'_{0}(s)}\sqrt{c}f^{-2}(\tau)(-D+\frac{c}{f^{2}(\tau)})^{-1/2}d\tau =\mid w_{0} \mid (= \mbox {length of}\; \gamma_{F_s} \mbox{for all $s$}), 
\end{equation}
where the integral is possibly considered under Convention \ref{cv1}. But the integrand and the right hand side in (\ref{e5}) are independent of $s$, thus,  $\tau'_{0}(s)=\pi_{I}\circ exp_{z_{0}}(\alpha_{i}(s))$ is constant, and
\begin{equation} \label{negra}
\frac{d}{ds}\mid_{s=0}\pi_{I}\circ exp_{z_{0}}(\alpha_{i}(s))=0.
\end{equation}
From (\ref{zz}) and (\ref{negra}) $v_{i}=\frac{d}{ds}\mid_{s=0}\alpha_{i}(s)$ yields a direction of conjugacy of $\gamma$, for any $i=1,\dots m'$, and it is clear from the construction that these $m'$ directions are independent. $\Box$

{\it Remark.} Note that the following case may hold: the point $x_0$ has a conjugate point $x_1$ along the $F-$ geodesic $\hat \gamma_{F}$, but if we consider any geodesic $\gamma$ emanating from $z_{0}=(\tau_{0},x_{0})$ which projects on $\hat{\gamma}_{F}$, the reparameterization $\gamma_{F}$ of $\hat{\gamma}_{F}$ does not  reach until $x_1$ and, so, there is no a conjugate point $z_1$ of $z_0$ along $\gamma$ which projects onto $x_1$. That is, the geodesic $\gamma$ ``escapes'' at the extremes of $I$ before $\hat \gamma_F$ reaches $x_1$. This possibility may happen, for example, when the spacetime is extendible through the extremes of $I$. But it does not necessarily happens because of this reason; in fact, de Sitter spacetime, where $f= \cosh$, is a simple counterexample (recall that if $\int_c^b f = \infty$ all null geodesics are future-complete \cite{Sa98} and the GRW spacetime not only is not extendible through $b$ as a $GRW$ spacetime but also it is not extendible as a spacetime; compare all this discussion with \cite[p. 73]{Uh}). When the fiber is weakly convex  the necessary and sufficient conditions to ensure that, for any geodesic $\gamma$ non-tangent to the base, $\gamma_{F}$ will cover all $\hat{\gamma}_{F}$ are the ``non--escape'' equalities  
\be
\label{pl}
\int^{c}_{a}f^{-2}(\frac{1}{f^{2}}+1)^{-1/2}=\infty \quad \int^{b}_{c}f^{-2}(\frac{1}{f^{2}}+1)^{-1/2}=\infty
\ee
for certain $c\in (a,b)$ (see \cite[Lemma 4]{FS}). Recall that this condition implies Condition (A) and, so, the spacetime will be geodesically connected.
Summing up:

\begin{corollary}\label{cUh}
 Consider a GRW spacetime with weakly convex fiber where the ``non--escape'' equalities (\ref{pl}) hold. Then the spacetime is geodesically connected and any causal geodesic $\gamma(t) = (\tau(t),\gamma_{F}(t))$ starting at $z_0$ have conjugate points in bijective correspondence (including multiplicities) with the conjugate points of the inextendible geodesic $\hat{\gamma}_{F}(r)$ obtained from the projection $\gamma_{F}(t)$ on the fiber. 
\end{corollary}

{\it Remark.} This result allows to extend, in our ambient, the ones by Uhlenbeck for null geodesics \cite{Uh} to all causal geodesics.
For instance, normalize all causal geodesics (non tangent to the base) such that $c\equiv(f^{4}\circ \tau)\cdot g(\frac{d\gamma_{F}}{dt},\frac{d\gamma_{F}}{dt})=1$ and choose $D\leq 0$; all future-pointing causal geodesics starting at $z_0 =(\tau_0, x_0)$ and having associated the fixed value of $D=g^{f}(\frac{d\gamma}{dt},\frac{d\gamma}{dt})$, are in bijective correspondence with the $F-$geodesics starting at $x_0$, being the conjugate points  preserved. So: 
\begin{quote}
{\em Under the assumptions of Corollary \ref{cUh}, and fixed $D\leq 0$, if $x_0$ and $x_1$ are not conjugate the loop space of $F$ is homotopic to a cell complex constructed with a cell for each causal $D$-geodesic (with $c=1$) from $z_0$ to the line $L_{x_1}=\{ (t,x_1): t\in I\}$ with the dimension of the cell equal to the index of the $D$-geodesic.}
\end{quote}
Recall that in \cite{Uh} the conformal invariance of null conjugate points is explicitly used, but this invariance does not hold for timelike geodesics (bidimensional anti-de Sitter spacetime, which is globally conformal to a strip in Lorentz-Minkowski spacetime, is a simple example); this makes necessary our approach.

Theorem \ref{t-3} and equalities (\ref{pl}) can be also combined to yield Morse relations as follows. Fix two non-conjugate points $z_0=(\tau_0,x_0), z'_0=(\tau'_0,x'_0)$ and a field ${\cal K}$. Let $\Omega(z_0,z'_0)$ (resp. 
$\Omega(x_0,x'_0)$) be the space of continuous paths joining $z_0, z'_0$ in $I\times F$ (resp. $x_0, x'_0$ in $F$). Let ${\cal P}_{z_0,z'_0}(t)$ (resp. 
${\cal P}_{x_0,x'_0}(t)$) be the Poincar\'e polinomial of $\Omega(z_0,z'_0)$ (resp. $\Omega(x_0,x'_0)$); that is, ${\cal P}_{z_0,z'_0}(t)$ is the formal series
$$ {\cal P}_{z_0,z'_0}(t)= \beta_0 + \beta_1 t + \beta_2 t^2 + \cdots $$
where $\beta_q$ is the $q$--th Betti number of $\Omega(z_0,z'_0)$ for homology with coefficients in ${\cal K}$, $\beta_q = dim H^q(\Omega(z_0,z'_0), {\cal K})$. Clearly, ${\cal P}_{z_0,z'_0}(t) \equiv 
{\cal P}_{x_0,x'_0}(t)$. Let 
$$ {\cal M}_{z_0,z'_0}(t)= \bar a_0 + \bar a_1 t + \bar a_2 t^2 + \cdots $$ 
$$ (\mbox{resp.} \; {\cal M}_{x_0,x'_0}(t)=  a_0 +  a_1 t +  a_2 t^2 + \cdots) $$
be the Morse polinomials of $z_0, z'_0$ (resp. $x_0, x'_0$), i.e. $\bar a_q$ (resp $a_q$) is the number of geodesics joining $z_0$ and $z'_0$ (resp. $x_0$ and $x'_0$) with Morse index equal to $q$, where the Morse index of a geodesic connecting two fixed non-conjugate points is the sum of the indexes of conjugate poins to the first point along the geodesic. Then, under the hypotheses of Theorem \ref{t2} and from Theorem \ref{t-3}:
\begin{equation} \label{st} 
a_q \leq \bar a_q + \bar a_{q+1}, \quad \forall q\geq 0, 
\end{equation}
\begin{equation} \label{stt} 
 \bar a_0 >0 \Rightarrow a_0 >0 ; \quad \bar a_q >0 \Rightarrow a_{q-1} + a_q >0 ; \quad \forall q \geq 1.
\end{equation}
In particular, if the polinomials are finite then  ${\cal M}_{z_0,z'_0}(t) \geq {\cal M}_{x_0, x'_0}(t), \forall t\geq 1$. But if $(F,g)$ is a complete Riemannian manifold, then the well-known Morse relations implies the existence of a formal polinomial, with non-negative integer coefficients ${\cal Q}(t)$ such that
\begin{equation} \label{ttt}
 {\cal M}_{x_0,x'_0}(t) = {\cal P}_{x_0,x'_0}(t) + (1+t) {\cal Q}(t). 
\end{equation}

{\it Remark.} In general, it is not true that $a_0 \geq \bar a_0$ or $a_{q-1} + a_q \geq \bar a_q$. Recall that many geodesics in the GRW spacetime connecting $z_0, z'_0$ may project on the same pregeodesic of $F$. A simple counterexample of this is de Sitter spacetime (with a straightforward modification, one can also get that hypotheses in Theorem \ref{t2} are fulfilled). So, inequalities (\ref{stt}) cannot be improved.

Summing up,
\begin{corollary} \label{cM}
In a globally hyperbolic GRW spacetime satisfying either Condition (A) or Condition (B) with $d_{a}$, $d_{b}$ (if defined) equal to infinity the Morse inequalities (\ref{st}), (\ref{stt})  (with (\ref{ttt})) hold. 

As a consequence, if the Morse polinomial ${\cal M}_{z_0,z'_0}(t)$ is finite then, for each pair of non--conjugate points $z_0 , z'_0$ there exist a  polinomial ${\cal Q}(t)$ with non--negative integer coefficients and computable from the fiber such that 
$$ {\cal M}_{z_0,z'_0}(t) \geq {\cal P}_{z_0,z'_0}(t) + (1+t) {\cal Q}(t), \quad \forall t \geq 1. $$

\end{corollary}

\section{Applications} \label{s6}
\subsection {Two dimensional case} \label{s6.1}

Next we will particularize previous results to bidimensional GRW spacetimes with strongly convex fiber (necessarily an interval $(J, dx^2)$). Recall that in this case the opposite metric $-g^f$ is also Lorentzian and, in fact, it corresponds to a static (standard) spacetime. The chronological relation  can be now extended for non-causally related points, just defining that two points are {\it spacelike} related if they are chronologically related for $-g^f$. In fact, we will simplify our terminology with the following 
(re-)definition.

\begin{definition}\label{d3} Consider a bidimensional GRW (or static) spacetime. 
Two points $(\tau_{0},x_{0}), (\tau'_{0},x'_{0}) $ are spacelike [resp. timelike, lightlike] related iff there exists a spacelike [resp. timelike, lightlike with non-vanishing derivative] 
curve  joining them.
\end{definition}

From a direct computation  (see also \cite[Th.3.3 and Lemma 3.5]{Sa97}) we have: 

\begin{lemma}\label{l4} Given $(\tau_{0},x_{0}), (\tau'_{0},x'_{0})\in (I\times J,-d\tau^{2}+f^{2}dx^2)$, they are

(i) spacelike related if and only if $\int_{\tau_{0}}^{\tau'_{0}}f^{-1}<d(x_{0},x'_{0})$ 

(ii) lightlike related if and only if $\int_{\tau_{0}}^{\tau'_{0}}f^{-1}=d(x_{0},x'_{0})$

(iii) timelike related if and only if $\int_{\tau_{0}}^{\tau'_{0}}f^{-1}>d(x_{0},x'_{0})$
\end{lemma}

Now, as a consequence of Lemma \ref{l4} and Theorem \ref{t1} we have:

\begin{corollary}\label{c2} In a GRW spacetime $(I\times J,-d\tau^{2}+f^{2}dx^2)$:

(i) If $(\tau_{0},x_{0})$, $(\tau'_{0},x'_{0})$ are timelike [resp lightlike] related then there exist a unique geodesic (necessarily timelike [resp lightlike]) which joins them.

(ii) All $(\tau_{0},x_{0})$, $(\tau'_{0},x'_{0})$ which are spacelike related can be joined by a geodesic (necessarily spacelike) if and only if Condition (C) or Condition (R) holds.
\end{corollary}

From Theorem \ref{t-3} and the fact that there are no conjugate points on a manifold of dimension 1, we have:

\begin{corollary} \label{c4} In a GRW spacetime $(I\times J,-d\tau^{2}+f^{2}dx^2)$ no geodesic $\gamma(t)= (\tau (t), \gamma_{F}(t))$ without zeroes in $d\tau /dt$  have  conjugate points. 

In particular, causal geodesics are free of conjugate points.
\end{corollary}

Now, consider a bidimensional static spacetime, say  $(K \times J \subseteq \R{}^2,  g_S= dy^2 - f^2(y) dx^2)$ where $g_S$ can be seen as the reversed metric of a GRW spacetime. Summarizing the conclusions of Lemma \ref{l4} and Corollaries \ref{c2}, \ref{c4}, the following extension of Theorem 1.1 in \cite{BGM} can be given (see also  \cite[Prop. 6.6]{GMPT}).

\begin{corollary}\label{c3} Given $(y_{0},x_{0}), (y'_{0},x'_{0})$ in the static spacetime $(K \times J \subseteq \R{}^2,  g_S= dy^2 - f^2(y) dx^2)$, they are

(i) spacelike related if and only if $\int_{y_{0}}^{y'_{0}}f^{-1}>d(x_{0},x'_{0})$. In this case there exist a unique  geodesic which joins them; this geodesic is necessarily spacelike and without conjugate points. 

(ii) lightlike related if and only if $\int_{y_{0}}^{y'_{0}}f^{-1}=d(x_{0},x'_{0})$. In this case there exist a unique  geodesic which joins them; this geodesic is necessarily  lightlike and without conjugate points. 

(iii) timelike related if and only if $\int_{y_{0}}^{y'_{0}}f^{-1}<d(x_{0},x'_{0})$. All points which are timelike related can be joined by a geodesic (necessarily timelike) if and only if Condition (C) or Condition (R) holds.
\end{corollary}

{\it Remark.} In fact,  no geodesic of the static spacetime without zeroes in the derivative of its spacelike component has conjugate points.
Anti de-Sitter spacetime is an example of static spacetime where all the   timelike geodesics have  conjugate points. Moreover, it is not geodesically connected.

\subsection{Conditions on  curvature} \label{s6.2}

As commented in the Introduction,  it is natural to assume, for a realistic GRW spacetime, that $Ric(\partial_t,\partial_t) \geq 0$, and it is straightforward to check that this condition is equivalent to $f''\leq 0$ (see \cite[Cor. 7.43]{O}). Recall that in this case $lim_{\tau\rightarrow a,b}f'$ and $lim_{\tau\rightarrow a,b}f$ always exist. So taking into account the cases in Table 1 we see that Condition (A) always holds except when $b<\infty$ (resp. $a>-\infty$) and $f'(b)>0$ (resp. $f'(a)<0$). In this case, although the GRW spacetime is not geodesically connected,  it is possible to extend the warping function $f$ through $b$ (resp. $a$) obtaining so a extended spacetime, which is also GRW. The GRW spacetime will be called {\it inextendible} if whenever an extreme of $I$ is finite, then $f$ cannot be extended continuously at these extremes to a real value $\alpha>0$. It seems clear that from a physical viewpoint just inextendible GRW spacetimes must be taken into account.

Therefore, Theorems \ref{t0} and \ref{t2} are appliable to these inextendible GRW spacetimes, yielding the points (ii) and (iv) in the following Corollary (the other two are included for the sake of completeness).

\begin{corollary}\label{c00} An inextendible GRW spacetime with $Ric(\partial_{t},\partial_{t})\geq 0$ and weakly convex fiber satisfies

(i) Each two causally related points can be joined with one non-spacelike geodesic, which is unique if the fiber is strongly convex.

(ii) The spacetime is geodesically connected. Moreover, each strip $(\hat a,\hat b)\times F\subset I\times F$, $a<\hat a<\hat b<b$ with the restricted metric is geodesically connected if and only if $f'(\hat a) \geq 0$ and $ f'(\hat b)\leq 0$ (i.e. $f'(\hat a)\cdot f'(\hat b)\leq 0$).

(iii) There exist a natural surjective map between geodesics connecting $z_{0}=(\tau_{0},x_{0})$,  $z'_{0}=(\tau'_{0},x'_{0})\in I\times F$ and F-geodesics connecting $x_{0}$ and $x'_{0}$. Under this map, when the geodesic connecting $z_{0}$ and $z'_{0}$ is causal then the multiplicity of its conjugate points is equal to the multiplicity for the corresponding geodesic connecting $x_{0}$, $ x'_{0}$.

(iv) If $(F,g)$ is complete and $F$ is not contractible in itself, then any $z_{0}, z'_{0}\in I\times F$ can be joined by means of infinitely many spacelike geodesics. If $x_{0}, x'_{0}$ are not conjugate there are at most finitely many causal geodesics connecting them.
\end{corollary}

For the last assertion (ii), recall that it is straightforward from Theorem \ref{t1} under strongly convexity. But, from the proof of this theorem, this assumption can be dropped because $f''\leq 0$ (recall that then Conditions (A), (B), (C) are equivalent and Condition (R) is not appliable).

Finally, we give a further consequence of equalities (\ref{pl}):

\begin{corollary} \label{c6}
Consider a GRW spacetime $(I\times F, g^f)$ which is globally hyperbolic  and satisfies the non-escape equalities (\ref{pl}), and fix $D_{0}\leq 0$. If any geodesic $\gamma(t)=(\tau(t),\gamma_{F}(t))$ starting at $z_0 =(t_0, x_0)$ and having associated values $D$, $c$ equal to $D_{0}$, $1$, respectively, is free of conjugate points then the fiber can be covered topologically by $\R{}^{n}$, being $n=dimF$. 
\end{corollary}

{\it Proof.} Under this assumption the F-geodesics starting at $x_{0}$ have no conjugate points and so, as $F$ is complete, exp$_{x_0}: T_{x_0}F \equiv \R{}^n \rightarrow F$ is a surjective local diffeomorphism. Taking the pull-back metric on $T_{x_0}F$, a local isometry with domain a complete manifold (and so a Riemannian covering) is obtained. $\Box$ 

{\it Remark.} The assumption on conjugate points when $D_{0}=0$ holds if in the future of $z_0$ we have $R(X,Y,Y,X) \leq 0$ whenever X, Y span a degenerate plane on a lightlike geodesic starting at $z_{0}$ (see \cite[Th. 10.77]{BEE}); moreover, the non-escape inequalities (\ref{pl}) can be reduced to

\be
\int^{c}_{a}f^{-1}=\infty \quad \int^{b}_{c}f^{-1}=\infty
\ee
when just null geodesics are considered, so we reobtain \cite[Theorem 5.3]{Uh} in our ambient.

\newpage
\begin{figure}
\begin{center}
\includegraphics[width=15cm]{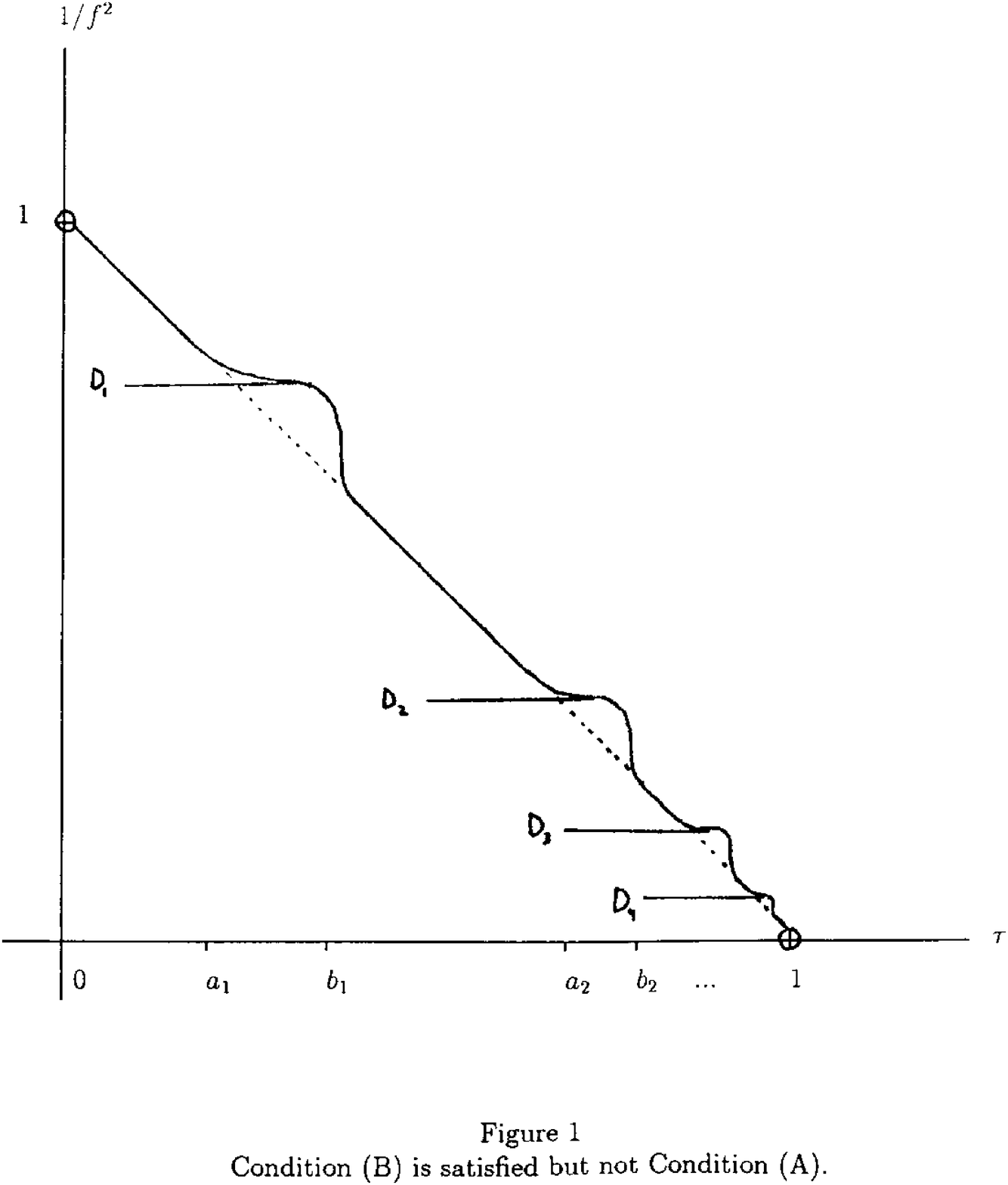}
\end{center}
\end{figure}

\newpage

\begin{figure}
\begin{center}
\includegraphics[width=10cm]{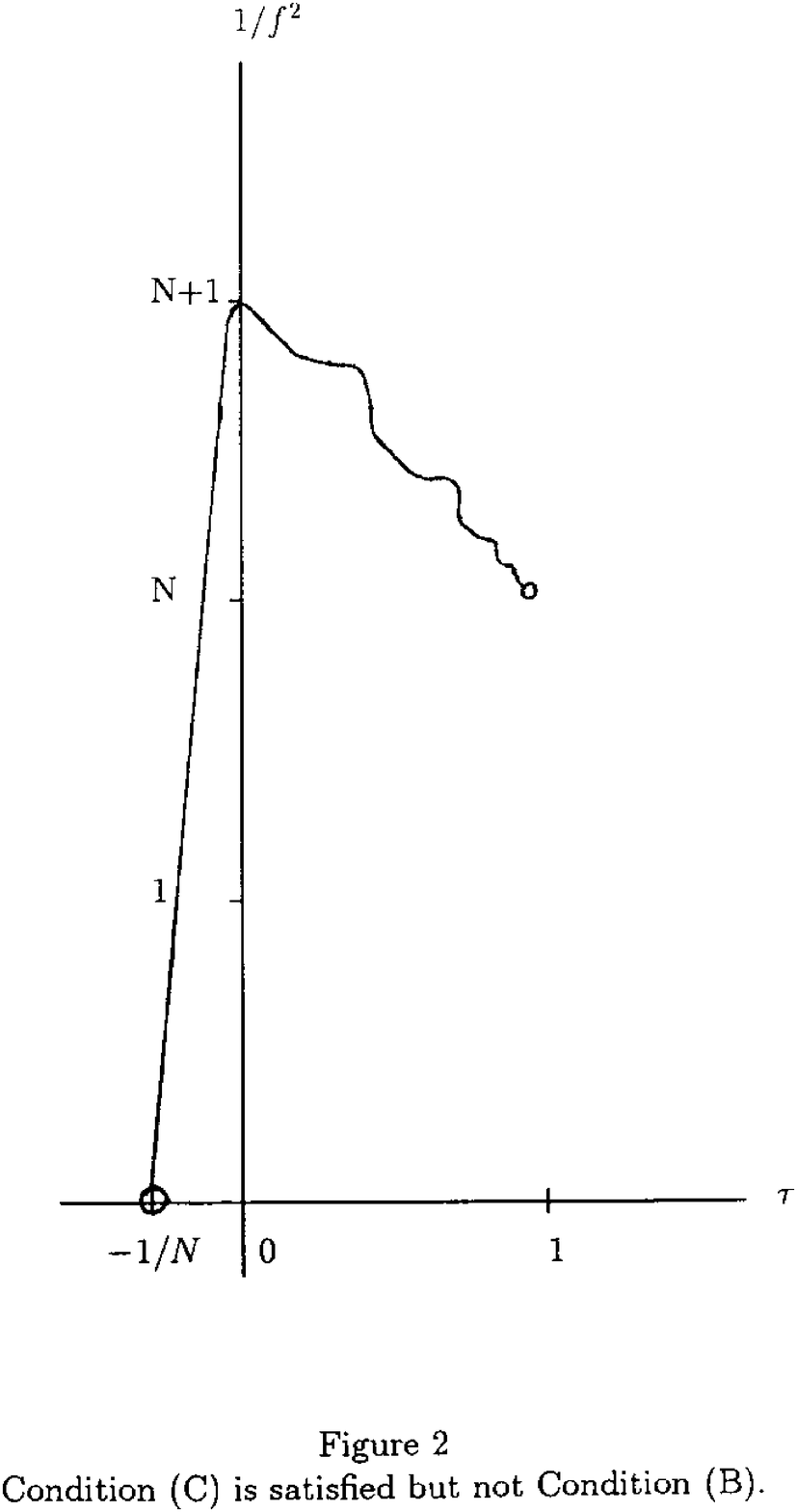}
\end{center}
\end{figure}

\newpage

\[%
\hspace*{-2.3cm}%
\begin{tabular}
[c]{c|cc|c|cc|c|}\cline{2-7}%
& \multicolumn{3}{|c}{b$<\infty$} & \multicolumn{3}{|c|}{b=$\infty$%
}\\\cline{2-7}\cline{2-7}%
& \multicolumn{2}{|c|}{$\lim\limits_{\tau\rightarrow b}f(\tau)$} &
\multicolumn{1}{|c|}{Condition (A)} & \multicolumn{2}{|c|}{$\lim
\limits_{\tau\rightarrow b}f(\tau)$} & Condition (A)\\\hline
\multicolumn{1}{|c|}{1} & \multicolumn{2}{|c|}{0} & \multicolumn{1}{|c|}{Yes}%
& \multicolumn{2}{|c|}{0} & Yes\\\hline
\multicolumn{1}{|c|}{} &  & \multicolumn{1}{|c|}{$f^{\prime}$ no extendible to
$b$} & \multicolumn{1}{|c|}{No information} &  &  & \\\cline{3-4}%
\multicolumn{1}{|c|}{} &  & \multicolumn{1}{|c|}{$\lim\limits_{\tau\rightarrow
b}f^{\prime}=\beta\in\lbrack-\infty,0)$} & \multicolumn{1}{|c|}{Yes} &  &  &
\\\cline{3-4}%
\multicolumn{1}{|c|}{2} & $\alpha\in\mathbb{R},\alpha\neq0$ &
\multicolumn{1}{|c|}{$\lim\limits_{\tau\rightarrow b}f^{\prime}=\beta
\in(0,\infty]$} & \multicolumn{1}{|c|}{No*} & \multicolumn{2}{|c|}{$\alpha
\in\mathbb{R},\alpha\neq0$} & Yes\\\cline{3-4}%
\multicolumn{1}{|c|}{} &  & \multicolumn{1}{|c|}{$\lim\limits_{\tau\rightarrow
b}f^{\prime}=0$ and $f^{\prime\prime}$} & \multicolumn{1}{|c|}{Yes} &  &  & \\
\multicolumn{1}{|c|}{} &  & \multicolumn{1}{|c|}{bounded in $[b-\epsilon,b)$}%
& \multicolumn{1}{|c|}{} &  &  & \\\hline
\multicolumn{1}{|c|}{3} & \multicolumn{2}{|c|}{$\infty$} &
\multicolumn{1}{|c|}{No} & $\infty$ & \multicolumn{1}{|c|}{$\int_{c}^{b}%
\frac{1}{f}=\infty$} & Yes\\\cline{6-7}%
\multicolumn{1}{|c|}{} &  &  & \multicolumn{1}{|c|}{} &  &
\multicolumn{1}{|c|}{$\int_{c}^{b}\frac{1}{f}<\infty$} & No**\\\hline
\end{tabular}
\]

\ \ \ \ \ \ \ \ \ \ \ \ \ \ \ \ \ \ \ \ \ \ \ \ \ \ \ \ \ \ \ \ \ \ \ \ \ \ \ Table
1 \ 

If $f$ is continuously extendible to $b$, when Condition (A) is satisfied at
$b$.

\bigskip

\bigskip

\bigskip(*) Condition (C) does not hold either. No information on Condition
(R), if appliable.

(**) No information on Condition (C) or (R).


\begin{thebibliography}{99}

\bibitem[ARS]{ARS} L.J. Al\'{\i}as, A. Romero, M. S\'{a}nchez, Uniqueness of complete spacelike hypersurfaces of constant mean curvature in Generalized Robertson-Walker spacetimes {\it Gen. Relat. Gravit.} {\bf 27} (1995) 71--84.

\bibitem[BEE]{BEE} J.K. Beem, P.E. Ehrlich and K.L. Easley, {\it Global Lorentzian geometry}, Monographs Textbooks Pure Appl. Math. {\bf 202} (Dekker Inc., New York, 1996).

\bibitem[BF]{BF} V. Benci, D. Fortunato,  Existence  of  geodesics 
for the Lorentz metric of a stationary gravitational field, {\it Ann. Inst. Henri Poincar\'e}, {\bf 7}  (1990) 27-35.

\bibitem[BGM]{BGM} V. Benci, F. Giannoni, A. Masiello, Some properties of the spectral flow in semiriemannian geometry, {\it J. Geom. Phys.} {\bf 27} (1998) 267--280. 

\bibitem[BGS]{BGS} R. Bartolo, A. Germinario and M. S\'anchez, Convexity of domains of Riemannian manifolds, {\it Ann. Global Anal. Geom.}, to appear.

\bibitem[BM]{BM} V. Benci,  A. Masiello, A Morse index for geodesics in static Lorentzian manifolds, {\it Math. Ann.} {\bf 293} (1992) 433--442. 

\bibitem[CM]{CM} E. Calabi, L. Markus, Relativistic space forms, {\it Ann. Math.} {\bf 75} (1962) 63-76.


\bibitem[FS]{FS} J.L. Flores, M. S\'anchez, Geodesic connectedness of multiwarped spacetimes, {\it J. Diff. Equat.}, to appear.

\bibitem[Gi]{Gi}  F.  Giannoni,  Geodesics  on  non  static  Lorentz  manifolds   of Reissner-Nordstr\"{o}m type, {\it Math. Ann.} {\bf 291} (1991) 383-401.

\bibitem[GM]{GM}  F. Giannoni, A. Massiello, Geodesics on  Lorentzian  manifolds  with quasi-convex boundary, {\it Manuscripta Math.} {\bf 78} (1993) 381-396.


\bibitem[GMPa]{GMP} F. Giannoni, A. Masiello, P. Piccione, Convexity and the finiteness of the number of geodesics: Applications to the multiple-image effect, {\it Class. Quant. Grav.} {\bf 16} (1999) 731--748.

\bibitem[GMPb]{GMPb} F. Giannoni, A. Masiello, P. Piccione, A Morse theory for massive particle and photons in General Relativity, {\it J. Geom. Phys.}, to appear.


\bibitem[GMPT]{GMPT} F. Giannoni, A. Masiello, P. Piccione, D. Tausk, A generalized index theorem for Morse-Sturm systems and applications to semi-Riemannian geometry, preprint (1999).

  \bibitem[He]{He} A.D. Helfer, Conjugate points on spacelike geodesics or pseudo-self-adjoint Morse-Sturm-Liouville systems, {\it Pac. J. Math.} {\bf 164} (1994) 321--350.

\bibitem[Ma]{Ma} A. Masiello, Variational methods in Lorentzian Geometry, Pitman Research Notes in Mathematics Series {\bf 309}, Longman Scientific and Technical, Harlow, Essex (1994).

\bibitem[O]{O} B. O'Neill, Semi-Riemannian Geometry with applications to Relativity, Series in Pure and Applied Math. {\bf 103} Academic Press, N.Y. (1983).

\bibitem[Sa97]{Sa97} M. S\'{a}nchez, Geodesic connectedness in Generalized Reissner-Nordstr\"{o}m type Lorentz manifolds, {\it Gen. Relat. Gravit.} {\bf 29} (1997) 1023--1037.

\bibitem[Sa98]{Sa98} M. S\'{a}nchez, On the Geometry of Generalized Robertson-Walker spacetimes: Geodesics, {\it Gen. Relat. Gravit.} {\bf 30} (1998) 915--932.


\bibitem[Sa99]{Sa99} M. S\'{a}nchez, On the Geometry of Generalized Robertson-Walker spacetimes: Curvature and Killing fields, {\it J. Geom. Phys.} {\bf 31} (1999) 1--15.

\bibitem[Sc]{Sc} H.-J. Schmidt, How should we measure spatial distances? {\it Gen. Relat. Gravit.} {\bf 28} No. 7 (1996) 899-903.

\bibitem[Uh]{Uh} K. Uhlenbeck, A Morse Theory for geodesics on Lorentz manifolds, {\it Topology} {\bf 14} (1975) 69--90. 

\end{thebibliography}
\end{document}